ANNALES
DE L'INSTITUT
HENRI
POINCARÉ
PROBABILITÉS
ET STATISTIQUES

# Moderate deviations for some point measures in geometric probability

## Yu Baryshnikov[a], P. Eichelsbacher[b], T. Schreiber[1,c] and J. E. Yukich[2,d]

[a]Rm 2C-323, Bell Laboratories, Lucent Technologies, 600-700 Mountain Ave, Murray Hill, NJ 07974, USA.
E-mail: ymb@research.bell-labs.com
[b]Fakultät für Mathematik, Ruhr-Universität Bochum, NA 3/68, 44780 Bochum, Germany.
E-mail: peter.eichelsbacher@rub.de
[c]Faculty of Mathematics and Computer Science, Nicholas Copernicus University, Toruń, Poland.
E-mail: tomeks@mat.uni.torun.pl
[d]Department of Mathematics, Lehigh University, Bethlehem, PA 18015, USA. E-mail: joseph.yukich@lehigh.edu

**Abstract.** Functionals in geometric probability are often expressed as sums of bounded functions exhibiting exponential stabilization. Methods based on cumulant techniques and exponential modifications of measures show that such functionals satisfy moderate deviation principles. This leads to moderate deviation principles and laws of the iterated logarithm for random packing models as well as for statistics associated with germ-grain models and $k$ nearest neighbor graphs.

**Résumé.** Les fonctionnelles en probabilité géométrique s'expriment souvent comme des sommes de fonctions bornées qui possèdent la fonction de stabilisation. Les méthodes de cumulants et les modifications exponentielles des mesures démontrent que ces fonctionnelles vérifient le principe des déviations modérées. Ceci donne des principes des déviations modérées et des lois de logarithme itéré pour les modèles de 'packing aléatoires' ainsi que pour des statistiques de modèles de 'germe-grain' et de graphes avec $k$ plus proches voisins.

## 1. Introduction

Functionals and measures induced by binomial and Poisson point processes in $d$-dimensional Euclidean space often satisfy a weak spatial dependence structure termed stabilization [27, 28], which, roughly speaking, quantifies the degree to which functionals are determined by the local configuration of points. Since the appearance of [27, 28], stabilization has been used in a general setting to establish thermodynamic limits [25, 29] and Gaussian limits for re-normalized functionals as well as re-normalized spatial point measures [4, 24, 26, 30]. Such general results can be applied to deduce limit laws for a variety of functionals and measures, including those defined by percolation models [24], random graphs in computational geometry [4, 24, 27], random packing models [3, 25, 28], germ-grain models [4, 25], as well as those involving maximal points [2] and vertices in convex hulls of i.i.d. samples [33].

---

[1]Research partially supported by the Foundation for Polish Science (FNP) and by the Polish Minister of Scientific Research and Information Technology, Grant 1 P03A 018 28 (2005-2007).
[2]The material is based upon work supported by the NSF under Grant DMS-0203720.

In this paper we use stabilization methods, cumulant techniques, and exponential modification of measures to establish asymptotics for random measures and functionals on scales intermediate between those appearing in Gaussian limit behavior and laws of large numbers. By appealing to Gärtner–Ellis and Dawson–Gärtner theory, this leads to moderate deviation principles and laws of the iterated logarithm for functionals of random sequential packing models as well as for statistics associated with germ-grain models and $k$ nearest neighbor graphs. By explicitly identifying rate functions we relate the large scale limit behavior of stabilizing functionals to the local behavior of the underlying density of points.

Recall that a family of probability measures $(\mu_\varepsilon)_{\varepsilon>0}$ on some topological space $\mathcal{T}$ obeys a large deviation principle (LDP) with speed $\varepsilon$ and good rate function $I(\cdot)\colon \mathcal{T}\to[0,\infty]$ if:

- $I$ is lower semi-continuous and has compact level sets $N_L := \{x\in\mathcal{T}\colon I(x)\le L\}$, for every $L\in[0,\infty)$,
- for every open set $G\subseteq\mathcal{T}$ we have

$$\liminf_{\varepsilon\to 0}\varepsilon\log\mu_\varepsilon(G) \ge -\inf_{x\in G} I(x), \tag{1.1}$$

- and for every closed set $A\subseteq\mathcal{T}$ we have

$$\limsup_{\varepsilon\to 0}\varepsilon\log\mu_\varepsilon(A) \le -\inf_{x\in A} I(x). \tag{1.2}$$

Similarly a family of random variables $(Y_\varepsilon)_{\varepsilon>0}$ with topological state space $\mathcal{T}$ obeys a LDP with speed $\varepsilon$ and good rate function $I(\cdot)\colon \mathcal{T}\to[0,\infty]$ if the sequence of their distributions obeys a LDP. We say that a sequence of random variables satisfies a moderate deviation principle (MDP) whenever the scaling is between that of an ordinary law of large numbers and that of a central limit theorem. Formally a moderate deviation principle is nothing but a LDP.

Let $\mathcal{C}([0,1]^d)$ be the collection of continuous $f\colon[0,1]^d\to\mathbb{R}$. For all $f\in\mathcal{C}([0,1]^d)$, $\|f\|_\infty$ denotes the essential supremum of $f$ and $\langle f,\mu\rangle$ denotes the integral of $f$ with respect to a signed finite variation Borel measure $\mu$ on $[0,1]^d$. For all $x\in\mathbb{R}^d$ and $r>0$, let $B_r(x)$ denote the Euclidean ball centered at $x$ of radius $r$. Denote the origin of $\mathbb{R}^d$ by $\mathbf{0}$. For all $\tau>0$, let $\mathcal{P}_\tau$ denote a homogeneous Poisson point process on $\mathbb{R}^d$ of intensity $\tau$. All random variables are defined on a common underlying probability space $(\Omega,\mathcal{F},\mathbf{P})$.

## 2. Random sequential packing

The following prototypical random sequential packing model arises in diverse disciplines, including physical, chemical, and biological processes. See [28] for a discussion of the many applications, the many references, and also a discussion of previous mathematical analysis. In one dimension, this model is often referred to as the Rényi car parking model [31].

With $N(\lambda)$ standing for a Poisson random variable with parameter $\lambda$, let $B_{\lambda,1}$, $B_{\lambda,2},\ldots,B_{\lambda,N(\lambda)}$ be a sequence of $d$-dimensional balls of volume $\lambda^{-1}$ whose centers are i.i.d. random $d$-vectors $X_1,\ldots,X_{N(\lambda)}$ with continuous probability density function $\kappa\colon[0,1]^d\to[0,\infty)$. Without loss of generality, assume that the balls are sequenced in the order determined by marks (time coordinates) in $[0,1]$. Let the first ball $B_{\lambda,1}$ be *packed*, and recursively for $i=2,3,\ldots,N(\lambda)$, let the $i$th ball $B_{\lambda,i}$ be packed iff $B_{\lambda,i}$ does not overlap any ball in $B_{\lambda,1},\ldots,B_{\lambda,i-1}$ which has already been packed. If not packed, the $i$th ball is discarded. The collection of centers of accepted balls induces a point measure on $[0,1]^d$. We call this the *random sequential packing measure* induced by balls (of volume $\lambda^{-1}$) with centers arising from $\kappa$.

For any finite point set $\mathcal{X}\subset\mathbb{R}^d$, assume the points $x\in\mathcal{X}$ have time coordinates which are independent and uniformly distributed over the interval $[0,1]$. Assume unit volume balls centered at the points of $\mathcal{X}$ arrive sequentially in an order determined by the time coordinates, and assume as before that each ball is packed or discarded according to whether or not it overlaps a previously packed ball. Let $\xi(x;\mathcal{X})$ be either 1 or 0 depending on whether the ball centered at $x$ is packed or discarded. Let $\xi_\lambda(x;\mathcal{X}) := \xi(\lambda^{1/d}x;\lambda^{1/d}\mathcal{X})$, where $\lambda^{1/d}x$ denotes scalar multiplication of $x$ and *not* the mark associated with $x$. Letting $\delta_x$ denote the Dirac point mass at $x$, the random sequential packing measures take the form

$$\mu_{\lambda\kappa}^\xi := \sum_{i=1}^{N(\lambda)} \xi_\lambda(X_i;\{X_j\}_{j=1}^{N(\lambda)})\delta_{X_i}. \tag{2.1}$$

Note that $\mu_{\lambda\kappa}^\xi$ is equal in distribution to $\sum_{x\in\mathcal{P}_{\lambda\kappa}} \xi_\lambda(x;\mathcal{P}_{\lambda\kappa})\delta_x$, where here and elsewhere $\mathcal{P}_{\lambda\kappa}$ denotes a Poisson point process on $[0,1]^d$ with intensity $\lambda\kappa$. For any random measure $\sigma$ on $\mathbb{R}^d$ we write $\bar{\sigma} := \sigma - \mathbb{E}[\sigma]$, so that for example $\bar{\mu}_{\lambda\kappa}^\xi := \mu_{\lambda\kappa}^\xi - \mathbb{E}[\mu_{\lambda\kappa}^\xi]$. Note that for all Borel sets $B \subset [0,1]^d$, we have $\mathbb{E}[\mu_{\lambda\kappa}^\xi(B)] := \lambda\int_B \mathbb{E}[\xi_\lambda(x;\mathcal{P}_{\lambda\kappa})]\kappa(x)\,dx$.

## 2.1. MDP for random sequential packing

From [4, 29], we know that $\xi$ depends upon local point configurations (formally termed "exponential stabilization" and defined in Section 6) and thus the one and two point correlation functions for $\xi_\lambda(x;\mathcal{P}_{\lambda\kappa})$ converge in the large $\lambda$ limit, which establishes volume order asymptotics for $\mathbb{E}[\mu_{\lambda\kappa}^\xi([0,1]^d)]$ and $\mathrm{Var}[\mu_{\lambda\kappa}^\xi([0,1]^d)]$ as $\lambda \to \infty$. Indeed, if for all $\tau > 0$, we put

$$V^\xi(\tau) := \mathbb{E}[\xi^2(\mathbf{0};\mathcal{P}_\tau)] + \int_{\mathbb{R}^d} (\mathbb{E}[\xi(\mathbf{0};\mathcal{P}_\tau \cup y) \cdot \xi(y;\mathcal{P}_\tau \cup \mathbf{0})] - \mathbb{E}[\xi(\mathbf{0};\mathcal{P}_\tau)]\mathbb{E}[\xi(y;\mathcal{P}_\tau)])\tau\,dy,$$

then [4, 26] we have

$$\lim_{\lambda\to\infty} \lambda^{-1}\mathrm{Var}[\mu_{\lambda\kappa}^\xi([0,1]^d)] = \int_{[0,1]^d} f^2(x)V^\xi(\kappa(x))\kappa(x)\,dx.$$

Additionally, the limit of the re-normalized measures $(\lambda^{-1/2}\bar{\mu}_{\lambda\kappa}^\xi)_\lambda$ is a generalized mean zero Gaussian field in the sense that the finite dimensional distributions of $(\lambda^{-1/2}\bar{\mu}_{\lambda\kappa}^\xi)_\lambda$ over test functions $f \in \mathcal{C}([0,1]^d)$ converge to those of a Gaussian field [4, 26].

It is natural to investigate the asymptotics of $(\bar{\mu}_{\lambda\kappa}^\xi)_\lambda$ on scales intermediate between those given by laws of large numbers and central limit theorems. Let $(\alpha_\lambda)_{\lambda>0}$ be such that $\lim_{\lambda\to\infty}\alpha_\lambda = \infty$ and $\lim_{\lambda\to\infty}\alpha_\lambda\lambda^{-1/2} = 0$.

We obtain the following MDP for packing measures:

**Theorem 2.1 (MDP on Poisson samples).** *For each $f \in \mathcal{C}([0,1]^d), f \not\equiv 0$, the family of random variables $(\alpha_\lambda^{-1}\lambda^{-1/2}\langle f,\bar{\mu}_{\lambda\kappa}^\xi\rangle)_\lambda$ satisfies on $\mathbb{R}$ the moderate deviation principle with speed $\alpha_\lambda^2$ and good rate function*

$$K_{\kappa;f}^\xi(t) := \frac{t^2}{2}\left(\int_{[0,1]^d} f^2(x)V^\xi(\kappa(x))\kappa(x)\,dx\right)^{-1}. \tag{2.2}$$

**Remarks.** (i) By taking $f \equiv 1$, Theorem 2.1 provides a MDP for the total number of balls accepted in the packing model with finite input. Theorem 2.1 adds to existing central limit theorems [3, 4, 8, 25, 28] and weak laws of large numbers [7, 28, 29] for random packing functionals. Through the rate function (2.2), Theorem 2.1 relates the large scale behavior of the family $(\alpha_\lambda^{-1}\lambda^{-1/2}\langle f,\bar{\mu}_{\lambda\kappa}^\xi\rangle)_\lambda$ to the local behavior of the underlying Poisson point process. By [4], Section 3.2, we have $\int_{[0,1]^d} V^\xi(\kappa(x))\kappa(x)\,dx > 0$ and thus the rate function is well-defined.

(ii) Our methods can be modified to show that Theorem 2.1 also holds whenever the support of $\kappa$ is a compact convex subset of $\mathbb{R}^d$ with non-empty interior.

(iii) We do not know how to prove the analog of Theorem 2.1 for de-Poissonized measures, that is to say when $N(\lambda)$ is replaced by $\lambda$.

The next result is a MDP on the level of *measures*. Denote by $\mathcal{M}([0,1]^d)$ the real vector space of finite variation signed measures on $[0,1]^d$. Equip $\mathcal{M}([0,1]^d)$ with the *weak topology* generated by the sets $\{U_{f,x,\delta}, f \in \mathcal{C}([0,1]^d), x \in \mathbb{R}, \delta > 0\}$, where

$$U_{f,x,\delta} := \{\nu \in \mathcal{M}([0,1]^d) \colon |\langle f, \nu \rangle - x| < \delta\}.$$

The Borel sigma field generated by the weak topology is denoted by $\mathcal{B}$. Since the collection of linear functionals $\{\nu \mapsto \langle f, \nu \rangle \colon f \in \mathcal{C}([0,1]^d)\}$ is separating in $\mathcal{M}([0,1]^d)$, it is well known that this topology makes $\mathcal{M}([0,1]^d)$ into a locally convex, Hausdorff topological vector space, whose topological dual is the preceding collection, hereafter identified with $\mathcal{C}([0,1]^d)$. In Section 2.2, in the context of the law of the iterated logarithm, we shall also endow $\mathcal{M}([0,1]^d)$ with an alternative weaker topology in order to make it metrizable.

**Theorem 2.2 (Measure level MDP).** *The family* $(\alpha_\lambda^{-1} \lambda^{-1/2} \bar{\mu}_{\lambda\kappa}^\xi)_\lambda$ *satisfies the moderate deviation principle on* $\mathcal{M}([0,1]^d)$, *endowed with the weak topology, with speed* $\alpha_\lambda^2$ *and a convex, good rate function*

$$I_\kappa^\xi(\nu) := \frac{1}{2} \int_{[0,1]^d} \left( \frac{\mathrm{d}\nu}{V^\xi(\kappa(x))\kappa(x)\,\mathrm{d}x} \right)^2 V^\xi(\kappa(x))\kappa(x)\,\mathrm{d}x, \tag{2.3}$$

*if* $\nu \in \mathcal{M}([0,1]^d)$ *is absolutely continuous with respect to* $V^\xi(\kappa(x))\kappa(x)\,\mathrm{d}x$, *and* $+\infty$ *otherwise.*

It is an easy observation to obtain a multi-dimensional version of Theorem 2.1.

**Theorem 2.3.** *For each linearly independent collection of continuous functions* $f_1, \ldots, f_l \colon [0,1]^d \to \mathbb{R}$, $l \in \mathbb{N}$, *the family*

$$(\alpha_\lambda^{-1} \lambda^{-1/2} (\langle f_1, \bar{\mu}_{\lambda\kappa}^\xi \rangle, \ldots, \langle f_l, \bar{\mu}_{\lambda\kappa}^\xi \rangle))_\lambda$$

*satisfies the moderate deviation principle on* $\mathbb{R}^l$ *with speed* $\alpha_\lambda^2$ *and a good rate function*

$$I_{\kappa, f_1, \ldots, f_l}^\xi(t) := \frac{1}{2} \langle t, C^{-1}(\xi, \kappa, f_1, \ldots, f_l) t \rangle, \quad t \in \mathbb{R}^l, \tag{2.4}$$

*where* $C(\xi, \kappa, f_1, \ldots, f_l)$ *denotes the covariance matrix with entries*

$$C_{ij}(\xi, \kappa, f_1, \ldots, f_l) := \int_{[0,1]^d} f_i(x) f_j(x) V^\xi(\kappa(x))\kappa(x)\,\mathrm{d}x.$$

Note that the linear independence of $f_1, \ldots, f_l$ guarantees that the matrix $C(\xi, \kappa, f_1, \ldots, f_l)$ is invertible so that $C^{-1}(\xi, \kappa, f_1, \ldots, f_l)$ is well defined.

**Remarks.** (i) *Starting with Theorem 2.3, we can alternatively apply Theorem 3.3 in [1] to get the measure-valued result, Theorem 2.2. See also [12] and [13], where this approach is applied to prove large and moderate deviations for empirical measures.*

(ii) *We expect that Theorem 2.2 holds with respect to the strong topology on* $\mathcal{M}([0,1]^d)$. *Proving this would necessitate showing that the upcoming Proposition 6.1 holds for all bounded functions on* $[0,1]^d$.

## 2.2. Laws of the iterated logarithm for random sequential packing

For all $\lambda > \mathrm{e}$, put $\alpha_\lambda := \sqrt{\log \log \lambda}$ and

$$\zeta_{\lambda\kappa}^\xi := \alpha_\lambda^{-1} \lambda^{-1/2} \bar{\mu}_{\lambda\kappa}^\xi. \tag{2.5}$$

Further, denote by $\mathcal{K}_\kappa^\xi$ the 'unit ball' for $I_\kappa^\xi$, given by

$$\mathcal{K}_\kappa^\xi := \{\nu \in \mathcal{M}([0,1]^d): I_\kappa^\xi(\nu) \le 1\}. \tag{2.6}$$

Throughout this section and the corresponding proofs we will endow $\mathcal{M}([0,1]^d)$ with a topology weaker than that introduced above so as to ensure metrizability, thus considerably simplifying the arguments. To this end, we consider a countable family $W := \{f_1, f_2, \ldots\}$ of continuous functions on $[0,1]^d$ uniformly dense in $\mathcal{C}([0,1]^d)$ and define the metric

$$\mathrm{dist}_W(\theta_1, \theta_2) := \sum_{k=1}^\infty \frac{1}{2^k \|f_k\|_\infty} \left| \int f_k \, \mathrm{d}\theta_1 - \int f_k \, \mathrm{d}\theta_2 \right|. \tag{2.7}$$

It is clear that $\mathrm{dist}_W(\cdot, \cdot)$ is a well-defined metric on $\mathcal{M}([0,1]^d)$, inducing topology which is weaker than the weak topology of Section 2.1. Moreover, $(\mathcal{M}([0,1]^d), \mathrm{dist}_W)$ is easily seen to be separable. Note also that *any* countable collection of functions in $\mathcal{C}([0,1]^d)$ can be extended to $W$ as above. We will show that Theorem 2.2 yields the following general Strassen-type law of the iterated logarithm (LIL), with all statements referring to the topology induced by $\mathrm{dist}_W(\cdot, \cdot)$.

**Theorem 2.4.** *For any possible coupling of the family of random measures $(\mu_{\lambda\kappa}^\xi)_\lambda$ on a common probability space $(\Omega, \mathcal{F}, \mathbf{P})$ and for any countable sequence $\lambda \to \infty$ the family of random measures $(\zeta_{\lambda\kappa}^\xi)_\lambda \subseteq \mathcal{M}([0,1]^d)$ is almost surely relatively compact and all its accumulation points almost surely fall into $\mathcal{K}_\kappa^\xi$. Moreover, there exists a coupling of $(\mu_{\lambda\kappa}^\xi)_\lambda$ on a common probability space such that the set of accumulation points of $(\zeta_{\lambda\kappa}^\xi)_\lambda$ almost surely coincides with $\mathcal{K}_\kappa^\xi$.*

It should be emphasized that we consider the families of random measures $(\zeta_{\lambda\kappa}^\xi)_\lambda$ along countable sequences $\lambda \to \infty$ rather than over all of $\mathbb{R}^+$ in order to avoid technicalities due to the presence of accumulation points arising along subsequences of $\lambda$ converging to a finite limit in $\mathbb{R}^+$. We do so in all of our LIL results below, without further mention.

The following scalar LIL is an immediate consequence of Theorem 2.4.

**Theorem 2.5 (LIL on Poisson samples).** *For any possible coupling of the family of random measures $(\mu_{\lambda\kappa}^\xi)_\lambda$ on a common probability space $(\Omega, \mathcal{F}, \mathbf{P})$ for any $f \in \mathcal{C}([0,1]^d)$ we have almost surely*

$$\limsup_{\lambda \to \infty} \langle f, \zeta_{\lambda\kappa}^\xi \rangle \le \sqrt{2 \int_{[0,1]^d} f^2(x) V^\xi(\kappa(x)) \kappa(x) \, \mathrm{d}x} \tag{2.8}$$

*and*

$$\liminf_{\lambda \to \infty} \langle f, \zeta_{\lambda\kappa}^\xi \rangle \ge -\sqrt{2 \int_{[0,1]^d} f^2(x) V^\xi(\kappa(x)) \kappa(x) \, \mathrm{d}x}. \tag{2.9}$$

*Moreover, there exists a coupling of $(\mu_{\lambda\kappa}^\xi)_\lambda$ on $(\Omega, \mathcal{F}, \mathbf{P})$ such that the above bounds are attained.*

"De-Poissonization" techniques for stabilizing functionals, as developed in [4, 27], yield a corresponding LIL for the measures generated by fixed-size binomial samples

$$\rho_{n,\kappa}^\xi := \sum_{i=1}^n \xi_n(X_i; \{X_j\}_{j=1}^n) \delta_{X_i}, \tag{2.10}$$

where $X_i$ are i.i.d. with density $\kappa$. Note that we are only able to state this result in the scalar setting. For notational convenience put

$$\theta_{n,\kappa}^\xi := \alpha_n^{-1} n^{-1/2} \bar{\rho}_{n,\kappa}^\xi. \tag{2.11}$$

For all locally finite $\mathcal{X} \subset \mathbb{R}^d$, let $H(\mathcal{X}) := H^\xi(\mathcal{X}) := \sum_{x \in \mathcal{X}} \xi(x; \mathcal{X})$. For all $\tau > 0$ consider the expected total effect of an inserted point at the origin on $H(\mathcal{P}_\tau)$. This is termed the mean "add one cost" and can be determined by either considering $\mathbb{E}[\Delta^\xi(\tau)]$ where

$$\Delta^\xi(\tau) := H(\mathcal{P}_\tau \cap B_S(\mathbf{0}) \cup \mathbf{0}) - H(\mathcal{P}_\tau \cap B_S(\mathbf{0})) \tag{2.12}$$

for $S$ large enough (see [4], Section 3.2) or by considering (see (2.16) of [26])

$$\delta^\xi(\tau) := \mathbb{E}[\xi(\mathbf{0}, \mathcal{P}_\tau)] + \tau \int_{\mathbb{R}^d} \mathbb{E}[\xi(\mathbf{0}, \mathcal{P}_\tau \cup y) - \xi(\mathbf{0}, \mathcal{P}_\tau)] \, \mathrm{d}y.$$

The mean add one cost is useful in evaluating $\mathrm{Var}[\langle f, \rho_{n,\kappa}^\xi \rangle]$: if for any $f \in \mathcal{C}([0,1]^d)$ we let

$$\sigma^2(\xi, \kappa, f) := \int_{[0,1]^d} f^2(x) V^\xi(\kappa(x)) \kappa(x) \, \mathrm{d}x - \left( \int_{[0,1]^d} f(x) \delta^\xi(\kappa(x)) \kappa(x) \, \mathrm{d}x \right)^2$$

then (see Theorem 3.4 of [4] and Theorem 6.2 of [26]) $\lim_{n \to \infty} n^{-1} \mathrm{Var}[\langle f, \rho_{n,\kappa}^\xi \rangle] = \sigma^2(\xi, \kappa, f)$. Our LIL goes as follows.

**Theorem 2.6 (LIL on binomial samples).** *For the sequence of random measures $(\rho_{n,\kappa}^\xi)_n$ and for any $f \in \mathcal{C}([0,1]^d)$ we have almost surely*

$$\limsup_{n \to \infty} \langle f, \theta_{n,\kappa}^\xi \rangle \leq \sqrt{2} \sigma(\xi, \kappa, f) \tag{2.13}$$

*and*

$$\liminf_{n \to \infty} \langle f, \theta_{n,\kappa}^\xi \rangle \geq -\sqrt{2} \sigma(\xi, \kappa, f). \tag{2.14}$$

**Remarks.** (i) *As evident from the expression for $\sigma(\xi, \kappa, f)$, Poissonization contributes extra randomness. As noted in Theorem 2.1 of [27] and Theorem 2.2 of [4], $\sigma(\xi, \kappa, f)$ is strictly positive.*

(ii) *It should be noted, as further discussed in the proof of Theorem 2.5 (see (7.5) there), that the coupling under which the bounds (2.8) and (2.9) are attained, in the special case of the uniform density $\kappa$, coincides with a certain natural coupling often appearing in applications.*

The proofs of Theorems 2.1–2.3 will be given in Sections 6.1 and 6.3, whereas the proofs of Theorems 2.4–2.6 are in Sections 7 and 8.

## 3. Spatial birth-growth models

Theorems 2.1–2.6 for the prototypical packing measures in Section 2 extend to measures arising from more general packing models. Consider for example the following spatial birth-growth model in $\mathbb{R}^d$. Let $\Psi := \{(X_i, T_i) \in \mathbb{R}^d \times [0,1]\}$ be a spatial-temporal Poisson point process. Seeds appear at random locations $X_i \in \mathbb{R}^d$ at times $T_i \in [0,1]$. When a seed is born at $X_i$ it has initial radius $\rho_i$, $0 < 1 \leq \rho_i \leq L < \infty$, and thereafter the radius grows at a constant speed $v_i$, generating a cell growing radially in all directions. When one expanding cell touches another, they both stop growing in their respective directions. In any event, we assume that the seed radii are deterministically bounded, i.e., they never exceed a fixed cut-off and they stop growing upon reaching it. Moreover, if a seed appears at $X_i$ and if the ball centered at $X_i$ with radius $\rho_i$ overlaps any of the existing cells, then the seed is discarded. Variants of this well-studied process are used to model crystal growth [34].

If seeds are born at random locations $X_i \in [0,1]^d$, it is natural to study the *spatial distribution* of accepted seeds. The convergence of the random measures induced by the locations of the accepted seeds is given by Theorem 3.5 of [4] and Theorem 2.1 of [28].

For any finite point set $\mathcal{X} \subset [0,1]^d$, assume the points $x \in \mathcal{X}$ have i.i.d. time marks over $[0,1]$. A mark at $x \in \mathcal{X}$ represents the arrival time of a seed at $x$. Assume that the seeds are centered at the points of $\mathcal{X}$, that they arrive sequentially in an order determined by the associated marks, and that each seed is accepted or rejected according to the rules above. Let $\xi(x; \mathcal{X})$ be either 1 or 0 according to whether the seed centered at $x$ is accepted or not. $\sum_{x \in \mathcal{X}} \xi(x; \mathcal{X})$ is the total number of seeds accepted.

As with the random sequential packing, let $X_1, \ldots, X_{N(\lambda)}$ be i.i.d. random variables with continuous density $\kappa$ on $[0,1]^d$ and with marks in $[0,1]$. The random measure

$$\mu_{\lambda\kappa}^{\xi} := \sum_{i=1}^{N(\lambda)} \xi_{\lambda}(X_i; \{X_j\}_{j=1}^{N(\lambda)}) \delta_{X_i}$$

is the scaled spatial birth-growth measure on $[0,1]^d$ induced by $X_1, \ldots, X_{N(\lambda)}$ and

$$\rho_{n,\kappa}^{\xi} := \sum_{i=1}^{n} \xi_n(X_i; \{X_j\}_{j=1}^{n}) \delta_{X_i}$$

is the scaled spatial birth-growth measure on $[0,1]^d$ induced by $X_1, \ldots, X_n$. Put

$$\zeta_{\lambda\kappa}^{\xi} := (\lambda \log\log \lambda)^{-1/2} \bar{\mu}_{\lambda\kappa}^{\xi}$$

and

$$\theta_{n,\kappa}^{\xi} := (n \log\log n)^{-1/2} \bar{\rho}_{n,\kappa}^{\xi}.$$

The following theorem is proved exactly along the lines of Theorems 2.1–2.6.

**Theorem 3.1.** *For each $f \in \mathcal{C}([0,1]^d)$, the family of random variables $(\alpha_{\lambda}^{-1} \lambda^{-1/2} \langle f, \bar{\mu}_{\lambda\kappa}^{\xi} \rangle)_{\lambda}$ satisfies the moderate deviation principle as in Theorem 2.1 whereas $\langle f, \zeta_{\lambda\kappa}^{\xi} \rangle_{\lambda}$ and $\langle f, \theta_{n,\kappa}^{\xi} \rangle_n$ satisfy the law of the iterated logarithm as in Theorems 2.5 and 2.6, respectively. Moreover, the family of measures $(\alpha_{\lambda}^{-1} \lambda^{-1/2} \bar{\mu}_{\lambda\kappa}^{\xi})_{\lambda}$ satisfies the moderate deviation principle on $\mathcal{M}([0,1]^d)$ with respect to the weak topology, as in Theorem 2.2 and the corresponding law of the iterated logarithm, as in Theorem 2.4.*

**Remarks.** (i) Theorem 3.1 adds to, [4, 5, 6, 28], which prove asymptotic normality for the number of accepted seeds.

(ii) Theorem 3.1 extends to more general versions of the prototypical packing model. The stabilization analysis of [28] yields MDPs and LILs in the finite input setting for the number of packed balls in the following general models: (a) models with balls replaced by particles of random size/shape/charge, (b) cooperative sequential adsorption models, and (c) ballistic deposition models (see [28] for a complete description of these models). In each case, our general MDP and LIL apply to the random packing measures associated with the centers of the packed balls, whenever the balls have a continuous density $\kappa : [0,1]^d \to [0,\infty)$.

## 4. Germ-grain models

Let $X_i$, $i \geq 1$, be i.i.d. with continuous density $\kappa$ on $[0,1]^d$. Let $T, T_i$, $i \geq 1$, be i.i.d. bounded random variables defined on the common probability space $(\Omega, \mathcal{F}, \mathbf{P})$, independent of the $X_i, i \geq 1$. Consider the random grains $X_i + n^{-1/d} B_{T_i}(\mathbf{0})$ as well as the random set

$$\varXi_n := \bigcup_{i=1}^{n} (X_i + n^{-1/d} B_{T_i}(\mathbf{0})),$$

where $B_r(x)$ again denotes the Euclidean ball centered at $x \in \mathbb{R}^d$ of radius $r > 0$. When the $X_i$, $1 \leq i \leq N(\lambda)$, are the realization of the Poisson point process $\mathcal{P}_{\lambda\kappa}$, the corresponding set $\Xi_{N(\lambda)}$ is a scale-changed Boolean model ([18], pp. 141, 233).

For all $u \in \mathbb{R}^d$, let $T(u)$ be a random variable with distribution equal to that of $T$. For all $x \in \mathbb{R}^d$ and all point sets $\mathcal{X} \subset \mathbb{R}^d$, denote by $V(x, \mathcal{X})$ the Voronoi cell around $x$ with respect to $\mathcal{X}$ and $L(x; \mathcal{X})$ the Lebesgue measure of the intersection of $\bigcup_{u \in \mathcal{X}} B_{T(u)}(u)$ and $V(x, \mathcal{X})$. Moreover, denote $L_\lambda(x; \mathcal{X}) := L(\lambda^{1/d}x; \lambda^{1/d}\mathcal{X})$.

The *volume measure* induced by $X_i$, $1 \leq i \leq N(\lambda)$, and $T_i$, $1 \leq i \leq N(\lambda)$, is

$$\mu_{\lambda\kappa}^L := \sum_{i=1}^{N(\lambda)} L_\lambda(X_i; \{X_j\}_{j=1}^{N(\lambda)})\delta_{X_i}.$$

Our next result gives a MDP for the volume measures. The proof is given in Section 6.2.

**Theorem 4.1.** *Assume that $\kappa$ is continuous and bounded away from zero on $[0,1]^d$. For each $f \in \mathcal{C}([0,1]^d)$, the family of random variables $(\alpha_\lambda^{-1}\lambda^{-1/2}\langle f, \bar{\mu}_{\lambda\kappa}^L\rangle)_\lambda$ satisfies the moderate deviation principle as in Theorem 2.1 whereas $\langle f, \zeta_{\lambda\kappa}^L\rangle_\lambda$ and $\langle f, \theta_{n,\kappa}^L\rangle_n$ satisfy the law of the iterated logarithm as in Theorems 2.5 and 2.6 respectively. Moreover, the family of measures $(\alpha_\lambda^{-1}\lambda^{-1/2}\bar{\mu}_{\lambda\kappa}^L)_\lambda$ satisfies the moderate deviation principle on $\mathcal{M}([0,1]^d)$ with respect to the weak topology, as in Theorem 2.2 and the corresponding law of the iterated logarithm as in Theorem 2.4.*

**Remark.** *Central limit theorems for volume measures are given by [4, 19, 25]. To the best of our knowledge there is no MDP result in the literature for the models considered here.*

## 5. $k$-nearest neighbors random graphs

Let $k$ be a positive integer. Given a locally finite point set $\mathcal{X} \subset \mathbb{R}^d$, the $k$-nearest neighbors (undirected) graph on $\mathcal{X}$, denoted $\mathrm{NG}(\mathcal{X})$, is the graph with vertex set $\mathcal{X}$ obtained by including $\{x, y\}$ as an edge whenever $y$ is one of the $k$-nearest neighbors of $x$ and/or $x$ is one of the $k$-nearest neighbors of $y$. The $k$-nearest neighbors (directed) graph on $\mathcal{X}$, denoted $NG'(\mathcal{X})$, is the graph with vertex set $\mathcal{X}$ obtained by placing a directed edge between each point and its $k$-nearest neighbors.

For all $t > 0$, let $\xi^t(x; \mathcal{X}) := 1$ if the length of the edge joining $x$ to its nearest neighbor in $\mathcal{X}$ is less than $t$ and zero otherwise. Put $\mu_{\lambda\kappa}^{\xi^t} := \sum_{x \in \mathcal{P}_{\lambda\kappa}} \xi_\lambda^t(x; \mathcal{P}_{\lambda\kappa})\delta_x$, $\zeta_{\lambda\kappa}^{\xi^t} := (\lambda \log \lambda)^{-1/2}\bar{\mu}_{\lambda\kappa}^{\xi^t}$, and $\theta_{n,\kappa}^{\xi^t} := (n \log \log n)^{-1/2}\bar{\rho}_{n,\kappa}^{\xi^t}$. We assume that $\kappa$ is continuous and bounded away from zero on $[0,1]^d$. The proof of the following MDP result is given in Section 6.3.

**Theorem 5.1.** *For each $f \in \mathcal{C}([0,1]^d)$ and each $t > 0$, the family of random variables $(\alpha_\lambda^{-1}\lambda^{-1/2}\langle f, \bar{\mu}_{\lambda\kappa}^{\xi^t}\rangle)_\lambda$ satisfies the moderate deviation principle as in Theorem 2.1 whereas $\langle f, \zeta_{\lambda\kappa}^{\xi^t}\rangle_\lambda$ and $\langle f, \theta_{n,\kappa}^{\xi^t}\rangle_n$ satisfy the law of the iterated logarithm as in Theorems 2.5 and 2.6, respectively. Moreover, the family of measures $(\alpha_\lambda^{-1}\lambda^{-1/2}\bar{\mu}_{\lambda\kappa}^{\xi^t})_\lambda$ satisfies the moderate deviation principle on $\mathcal{M}([0,1]^d)$ with respect to the weak topology, as in Theorem 2.2 and the corresponding law of the iterated logarithm as in Theorem 2.4.*

**Remarks.** (i) Theorem 5.1 yields a MDP and LIL for the empirical distribution function of the rescaled lengths of the edges in the nearest neighbors graph on $\mathcal{P}_{\lambda\kappa}$. This gives a MDP and LIL for the number of pairs of rescaled points distant at most $t$ from each other, adding to central limit theorems of (Chapter 4 of [23]).

(ii) Alternatively, given $m \in \mathbb{N}$, we could let $\xi_m^{NG}(x; \mathcal{X})$ (respectively, $\xi_m^{NG'}(x; \mathcal{X})$) be one or zero according to whether the degree of $x$ is equal to $m$. Then Theorem 5.1 holds for this definition of $\xi$, yielding a MDP and LIL for the number of vertices in the $k$-nearest neighbor graph of fixed degree.

## 6. Proof of the moderate deviations principles

### 6.1. Proof of the moderate deviations principles on Poisson samples

The proof of Theorem 2.1 as well as the Poissonized versions of Theorems 3.1, 4.1 and 5.1 involves the logarithmic Laplace transform of $\alpha_\lambda \lambda^{-1/2} \langle f, \bar{\mu}^\xi_{\lambda\kappa} \rangle, f \in \mathcal{C}([0,1]^d)$, defined by

$$\Lambda^\xi_{\lambda\kappa,\alpha_\lambda}(f) := \alpha_\lambda^{-2} \log[\mathbb{E} \exp(\alpha_\lambda \lambda^{-1/2} \langle f, \bar{\mu}^\xi_{\lambda\kappa} \rangle)]. \tag{6.1}$$

The existence of the Laplace transform follows easily from the boundedness of the $\xi$-functional (in the case of random sequential packing and its variants as well as for the $k$-nearest points graphs) or from the volume-order bound for the total mass of $\mu^\xi_{\lambda\kappa}$ (for germ-grain models). To prove Theorem 2.1 as well as the Poisson sample moderate deviations principles in Theorems 3.1, 4.1 and 5.1 it will suffice to establish the following result in each model.

**Proposition 6.1.** *For all models in Sections 2–5 and for all $f \in \mathcal{C}([0,1]^d)$ the logarithmic Laplace transform $\Lambda^\xi_{\lambda\kappa,\alpha_\lambda}(\cdot)$ satisfies*

$$\Lambda^\xi_\kappa(f) := \lim_{\lambda \to \infty} \Lambda^\xi_{\lambda\kappa,\alpha_\lambda}(f) = \frac{1}{2} \int_{[0,1]^d} f^2(x) V^\xi(\kappa(x)) \kappa(x) \, \mathrm{d}x.$$

**Proof of the Poisson sample moderate deviations principles.** Whenever the logarithmic Laplace transform of $\alpha_\lambda \lambda^{-1/2} \langle f, \bar{\mu}^\xi_{\lambda\kappa} \rangle$ satisfies the asymptotics of Proposition 6.1, the moderate deviations of the family $(\alpha_\lambda^{-1} \lambda^{-1/2} \langle f, \bar{\mu}^\xi_{\lambda\kappa} \rangle)_\lambda$ can be obtained via the following standard arguments. Fix $f \in \mathcal{C}([0,1]^d)$ and use Proposition 6.1 to get for all $s \in \mathbb{R}$

$$\Lambda^\xi_\kappa(sf) := \lim_{\lambda \to \infty} \Lambda^\xi_{\lambda\kappa,\alpha_\lambda}(sf) = \frac{s^2}{2} \int_{[0,1]^d} f^2(x) V^\xi(\kappa(x)) \kappa(x) \, \mathrm{d}x. \tag{6.2}$$

In particular, $\Lambda^\xi_\kappa(sf)$ is finite for all $s \in \mathbb{R}$ and, moreover, it is everywhere differentiable. Therefore, by the standard Gärtner–Ellis result (cf. Theorem 2.3.6 in [10]), the family $(\alpha_\lambda^{-1} \lambda^{-1/2} \langle f, \bar{\mu}^\xi_{\lambda\kappa} \rangle)_\lambda$ satisfies on $\mathbb{R}$ the full moderate deviation principle with speed $\alpha_\lambda^2$ and good rate function (Legendre–Fenchel-transform of $\Lambda^\xi_\kappa$)

$$K^\xi_{\kappa,f}(t) := \sup_{s \in \mathbb{R}}(ts - \Lambda^\xi_\kappa(sf)) = \frac{t^2}{2} \left( \int_{[0,1]^d} f^2(x) V^\xi(\kappa(x)) \kappa(x) \, \mathrm{d}x \right)^{-1}$$

as in (2.2). This yields Theorem 2.1 and the Poisson sample moderate deviations principle in Theorems 3.1, 4.1 and 5.1. □

It remains therefore to prove Proposition 6.1. To this end we shall consider appropriate exponential (Gibbsian) modifications of the considered point process $\mathcal{P}_{\lambda\kappa}$, with the energy functional given by the (negative) integral of a continuous function on $[0,1]^d$ against the empirical measure $\mu^\xi_{\lambda\kappa}$. In the context of packing functionals we will argue that the exponential stabilization property of the considered packing functionals can be used, thus enabling us to conclude the assertion of Proposition 6.1 using the method of cumulants and cluster measures developed in [4] in the context of the central limit theorem. Separate arguments will be needed for the models in sections four and five. Let us recall the exponential stabilization property [4].

### Stabilization

Exponential stabilization, used heavily in [4, 26, 27, 30], plays a central role in all that follows. For all $\tau > 0$, recall that $\mathcal{P}_\tau$ denotes a homogeneous Poisson point process on $\mathbb{R}^d$ of intensity $\tau$. The following is a slightly strengthened version of stabilization used in [4].

**Definition 6.1 ([4]).** *The functional $\xi$ is exponentially stabilizing between intensities $a$ and $b$, $0 < a < b < \infty$, if for all $x \in \lambda^{1/d}[0,1]^d$ and all $\lambda > 0$ there exists a random variable $R := R^{\xi}_{a,b,\lambda}(x)$ (a radius of stabilization for $\xi$ at $x$ between intensities $a$ and $b$) such that*

$$\xi(x;(\hat{\mathcal{P}} \cap B_R(x)) \cup \mathcal{X}) = \xi(x;(\hat{\mathcal{P}} \cap B_R(x)))$$

*for all finite $\mathcal{X} \subset \mathbb{R}^d \setminus B_R(x)$ and for all $\mathcal{P}_a \subseteq \hat{\mathcal{P}} \subseteq \mathcal{P}_b$ and, moreover, there exist finite constants $L := L(a,b) > 0$, $\alpha := \alpha(a,b) > 0$ such that for all $t > 0$*

$$\sup_{x \in \mathbb{R}^d, \lambda > 0} \mathbf{P}[R^{\xi}_{a,b,\lambda}(x) > t] \leq L \exp(-\alpha t). \tag{6.3}$$

*When $\xi$ stabilizes then for all $\tau > 0$ we let $\xi(\mathbf{0}; \mathcal{P}_\tau) := \lim_{l \to \infty} \xi(\mathbf{0}; \mathcal{P}_\tau \cap B_l(\mathbf{0}))$.*

Thus $R := R^{\xi}_{a,b,\lambda}(x)$ is a radius of stabilization if the value of $\xi(x; \hat{\mathcal{P}})$, with $\mathcal{P}_a \subseteq \hat{\mathcal{P}} \subseteq \mathcal{P}_b$, is unaffected by changes outside $B_R(x)$. For all classes of models considered here, the corresponding functionals $\xi$ are exponentially stabilizing, see [4].

*Moment and cumulant measures*

To put the ideas in precise terms, we first recall the formal definition of cumulants in the context specified for our purposes. Take $f \in \mathcal{C}([0,1]^d)$ and denote $W := [0,1]^d$ for short. Expand the *Laplace transform* $\mathbb{E} \exp(\alpha_\lambda \lambda^{-1/2} \langle f, \bar{\mu}^{\xi}_{\lambda\kappa} \rangle)$ in a power series in $f$ as follows:

$$\mathbb{E} \exp(\alpha_\lambda \lambda^{-1/2} \langle f, \bar{\mu}^{\xi}_{\lambda\kappa} \rangle) = 1 + \sum_{k=1}^{\infty} \frac{(\alpha_\lambda \lambda^{-1/2})^k \langle f^{\otimes k}, M^k_\lambda \rangle}{k!}, \tag{6.4}$$

where $f^{\otimes k} : \mathbb{R}^{dk} \to \mathbb{R}$, $k = 1, 2, \ldots$, is given by $f^{\otimes k}(v_1, \ldots, v_k) = f(v_1) \cdots f(v_k)$, and $v_i \in W, 1 \leq i \leq k$. $M^k_\lambda := M^k_\lambda(\kappa)$ is a measure on $\mathbb{R}^{dk}$, *the $k$th moment measure* (p. 130 of [9]). Both the existence of the moment measures and the convergence of the series (6.4) are direct consequences of the boundedness of $\xi$, given in all the models considered here.

Expanding the logarithm of the Laplace transform in the power series gives

$$\log \left[ 1 + \sum_{k=1}^{\infty} \frac{(\alpha_\lambda \lambda^{-1/2})^k \langle f^{\otimes k}, M^k_\lambda \rangle}{k!} \right] = \sum_{l=1}^{\infty} \frac{(\alpha_\lambda \lambda^{-1/2})^l \langle f^{\otimes l}, c^l_\lambda \rangle}{l!}; \tag{6.5}$$

the signed measures $c^l_\lambda$ are *cumulant measures* [22].

Note that the first cumulant measure coincides with the expectation measure and the second cumulant measure coincides with the covariance measure.

*Exponentially modified point processes*

For each $f \in \mathcal{C}([0,1]^d)$ we consider the Gibbs-modified point process $\mathcal{P}^{f \circ \xi}_{\lambda\kappa}$ given in law by

$$\frac{\mathrm{d}\mathcal{L}(\mathcal{P}^{f \circ \xi}_{\lambda\kappa})}{\mathrm{d}\mathcal{L}(\mathcal{P}_{\lambda\kappa})}[\mathcal{X}] := \frac{\exp(\sum_{x \in \mathcal{X}} f(x)\xi_\lambda(x, \mathcal{X}))}{\mathbb{E}\exp(\sum_{x \in \mathcal{P}_{\lambda\kappa}} f(x)\xi_\lambda(x, \mathcal{P}_{\lambda\kappa}))} = \frac{\exp(\sum_{x \in \mathcal{X}} f(x)\xi_\lambda(x, \mathcal{X}))}{\mathbb{E}\exp(\langle f, \mu^{\xi}_{\lambda\kappa} \rangle)} \tag{6.6}$$

for all finite point configurations $\mathcal{X} \subseteq [0,1]^d$, with $\mathcal{L}(\cdot)$ standing for the distribution of the argument random object.

Recall that if $Y$ is a real-valued random variable and $L(h) := \log \mathbb{E} \exp(hY)$ its log-Laplace transform, and if $(Y_h)_{h>0}$ are random variables with law

$$\mathbf{P}[Y_h \in \mathrm{d}y] := \frac{\exp(hy)}{\exp(L(h))} \mathbf{P}[Y \in \mathrm{d}y],$$

then the $k$th derivative of $L(h)$ at $h$ coincides with the $k$th cumulant of $Y_h$. Considering $Y := \langle f, \mu_{\lambda\kappa}^\xi \rangle$ and the Gibbs-modified point process $\mathcal{P}_{\lambda\kappa}^{hf\circ\xi}$ given by (6.6), we thus make the simple yet crucial observation that

$$\frac{\partial}{\partial h} \log \mathbb{E} \exp(h\langle f, \mu_{\lambda\kappa}^\xi \rangle) = \mathbb{E} \sum_{x \in \mathcal{P}_{\lambda\kappa}^{hf\circ\xi}} hf(x)\xi_\lambda(x, \mathcal{P}_{\lambda\kappa}^{hf\circ\xi})$$

and hence

$$\frac{\partial}{\partial h} \log \mathbb{E} \exp(h\langle f, \bar{\mu}_{\lambda\kappa}^\xi \rangle)_{|h=0} = 0. \tag{6.7}$$

Moreover,

$$\frac{\partial^2}{\partial h^2} \log \mathbb{E} \exp(h\langle f, \bar{\mu}_{\lambda\kappa}^\xi \rangle) = \frac{\partial^2}{\partial h^2} \log \mathbb{E} \exp(h\langle f, \mu_{\lambda\kappa}^\xi \rangle) = \mathrm{Var} \left[ \sum_{x \in \mathcal{P}_{\lambda\kappa}^{hf\circ\xi}} hf(x)\xi_\lambda(x, \mathcal{P}_{\lambda\kappa}^{hf\circ\xi}) \right]. \tag{6.8}$$

More generally, for all $k \geq 2$

$$\frac{\partial^k}{\partial h^k} \log \mathbb{E} \exp(h\langle f, \bar{\mu}_{\lambda\kappa}^\xi \rangle) = \frac{\partial^k}{\partial h^k} \log \mathbb{E} \exp(h\langle f, \mu_{\lambda\kappa}^\xi \rangle) = \langle f^{\otimes k}, c_{\lambda;hf}^k \rangle, \tag{6.9}$$

where $c_{\lambda;hf}^k$ is the $k$th cumulant measure of the empirical measure

$$\mu_{\lambda\kappa}^{hf\circ\xi} := \sum_{x \in \mathcal{P}_{\lambda\kappa}^{hf\circ\xi}} \xi_\lambda(x, \mathcal{P}_{\lambda\kappa}^{hf\circ\xi}) \tag{6.10}$$

or, for $k \geq 2$, equivalently of its centered version $\bar{\mu}_{\lambda\kappa}^{hf\circ\xi} := \mu_{\lambda\kappa}^{hf\circ\xi} - \mathbb{E}\mu_{\lambda\kappa}^{hf\circ\xi}$. The above observations (6.7) and (6.8) and the relation (6.9) with $h := \alpha_\lambda \lambda^{-1/2}$ leads to the following Taylor expansion of the log-Laplace transform $\Lambda_{\lambda\kappa,\alpha_\lambda}^\xi(tf)$ around zero, evaluated at $t = 1$:

$$\Lambda_{\lambda\kappa,\alpha_\lambda}^\xi(f) = \frac{1}{\alpha_\lambda^2} \left[ \frac{\alpha_\lambda^2}{2\lambda} \langle f^{\otimes 2}, c_\lambda^2 \rangle + \frac{\alpha_\lambda^3}{6\lambda^{3/2}} \langle f^{\otimes 3}, c_{\lambda;hf}^3 \rangle \right] = \frac{1}{2\lambda} \langle f^{\otimes 2}, c_\lambda^2 \rangle + \frac{\alpha_\lambda}{6\lambda^{3/2}} \langle f^{\otimes 3}, c_{\lambda;hf}^3 \rangle \tag{6.11}$$

for some $\eta \in [0,1]$. Note that for $\lambda$ large enough, $h = \alpha_\lambda \lambda^{-1/2}$ eventually falls into $[0,1]$. Below, we shall write

$$u := u(\eta, \lambda) := \eta h = \eta \alpha_\lambda \lambda^{-1/2}.$$

Now, from [4, 26, 28] all models in Sections 2–4 satisfy

$$\lim_{\lambda \to \infty} \lambda^{-1} \langle f^{\otimes 2}, c_\lambda^2 \rangle = \lim_{\lambda \to \infty} \lambda^{-1} \mathrm{Var}[\langle f, \bar{\mu}_\lambda^\xi \rangle] = \int_{[0,1]^d} f^2(x) V^\xi(\kappa(x))\kappa(x)\,dx.$$

Thus the second-order term on the right-hand side of (6.11) satisfies

$$\lim_{\lambda \to \infty} \frac{1}{2\lambda} \langle f^{\otimes 2}, c_\lambda^2 \rangle = \frac{1}{2} \int_{[0,1]^d} f^2(x) V^\xi(\kappa(x))\kappa(x)\,dx.$$

Therefore, in order to establish Proposition 6.1 it remains to show that the third-order term is negligible, i.e.

$$\lim_{\lambda \to \infty} \frac{\alpha_\lambda}{6\lambda^{3/2}} \langle f^{\otimes 3}, c_{\lambda;\eta hf}^3 \rangle = 0. \tag{6.12}$$

### 6.2. *Proof of* (6.12) *for random sequential packing and spatial birth-growth models*

In order to avoid unnecessary technicalities, assume without loss of generality that $\lambda = (Lm)^d$ for some large fixed $L > 4$ and some $m \in \mathbb{N}$. We partition the cube $[0,1]^d$ into $m^d := \lambda/L^d$ equal-sized sub-cubes $Q_{\bar{i}}^{(\lambda)}$, $\bar{i} \in \{0, \ldots, m\}^d$, with $Q_{\bar{i}}^{(\lambda)}$, $\bar{i} := (i_1, \ldots, i_d)$, arising as the translate of $[0, m^{-1}]^d = [0, \lambda^{-1/d}L]^d$ by $(Li_1, \ldots, Li_d)$. We shall call $\bar{i} := (i_1, \ldots, i_d)$ and $\bar{j} := (j_1, \ldots, j_d)$ neighbors iff $\max_{l=1}^{d} |i_l - j_l| \leq 1$, writing $\bar{i} \sim \bar{j}$. To establish the required relation (6.12), we require the following lemma, allowing us to control the difference between the Gibbs-modified process $\mathcal{P}_{\lambda\kappa}^{ufo\xi}$ with $u := u(\eta, \lambda) := \eta\alpha_\lambda\lambda^{-1/2}$ and the original process $\mathcal{P}_{\lambda\kappa}$.

**Lemma 6.1 (Coupling).** *For each $\delta > 0$, with $L$ chosen large enough and for each sufficiently large $\lambda$ there exists a probability space carrying versions of $\mathcal{P}_{\lambda\kappa}$ and $\mathcal{P}_{\lambda\kappa}^{ufo\xi}$ as well as a family $(\pi_{\bar{i}})_{\bar{i}\in\mathbb{Z}^d}$ of i.i.d. $\{0,1\}$-valued random variables with $\mathbf{P}[\pi_{\bar{i}} = 1] = \delta$ such that, with probability 1, whenever $\pi_{\bar{i}} = 0$ for some $\bar{i} \in \{0, \ldots, m\}^d$, we have the point processes $\mathcal{P}_{\lambda\kappa}^{ufo\xi}$ and $\mathcal{P}_{\lambda\kappa}$ coinciding on the corresponding sub-cube $Q_{\bar{i}}$ both in point locations and in mark values (given by the corresponding values of the re-scaled packing functional $\xi_\lambda$).*

**Proof.** We will apply the *domination by product measures result* of [21], more precisely Theorem 0.0 in [21]. It tells us that, for a family of $\{0,1\}$-valued random variables indexed by lattice vertices, if we are able to show that for each given site the probability of seeing 1 there conditioned on the configuration outside a fixed size neighborhood of the site exceeds certain large enough $p$, then this random field dominates a product measure with positive density $q$ which can be made arbitrarily close to 1 by appropriate choice of $p$. Throughout this proof and in the conclusion of the proof of (6.12) below, *whenever referring to point processes $\mathcal{P}_{\lambda\kappa}$ and $\mathcal{P}_{\lambda\kappa}^{ufo\xi}$ we have in mind their versions marked with their arrival times and with the corresponding values (zero or one) of the re-scaled packing functional $\xi_\lambda$.* In particular, the equality of such processes over a subdomain of $[0,1]^d$ is equivalent to the equality of point locations and arrival times and of their respective values of $\xi_\lambda$.

By the argument in Section 4 of [28] establishing the exponential stabilization of $\xi$, it follows there exists $L$ large enough such that for each $\bar{i} \in \{0, \ldots, m\}^d$ and for each marked point configuration $\eta$ outside $\bigcup_{\bar{j}\sim\bar{i}} Q_{\bar{j}}^{(\lambda)}$ we have that the variational distance between

(i) the law of the marked point process $\mathcal{P}_{\lambda\kappa}$ restricted to $Q_{\bar{i}}^{(\lambda)}$ and

(ii) the corresponding restriction of its conditional distribution given $\eta$

is uniformly bounded above by $L^d \exp(-\Omega(L))$. Indeed, this comes directly from bounding above the probability that either the stabilization sphere of at least one point of $\lambda^{1/d}(\mathcal{P}_{\lambda\kappa} \cap Q_{\bar{i}}^{(\lambda)})$ reaches $[\bigcup_{\bar{j}\sim\bar{i}} \lambda^{1/d} Q_{\bar{j}}^{(\lambda)}]^c$ or that a causal chain (i.e. a sequence of balls with increasing arrival times and with each successor overlapping its predecessor, see Section 4 of [28] for terminology) joining $\lambda^{1/d} Q^{(\lambda)}$ to $[\bigcup_{\bar{j}\sim\bar{i}} \lambda^{1/d} Q_{\bar{j}}^{(\lambda)}]^c$ arises due to the choice of $\eta$, see ibidem.

Given $\varepsilon > 0$ there is a $L := L(\varepsilon)$ large enough to guarantee the above variational distance is smaller than $\varepsilon/2$.

The next observation is that, with $u$ small enough, the conditional law of the Gibbs-modified marked point process $\mathcal{P}_{\lambda\kappa}^{ufo\xi}$ restricted to $Q_{\bar{i}}^{(\lambda)}$ given a marked point configuration $\eta$ outside $\bigcup_{\bar{j}\sim\bar{i}} Q_{\bar{j}}^{(\lambda)}$ differs in total variation from the corresponding conditional law for $\mathcal{P}_{\lambda\kappa}$ by at most $\varepsilon/2$. This is due to the fact that, by the definition (6.6), the density (Radon–Nikodym derivative) of the former conditional law with respect to the latter one at a point configuration $\mathcal{X} \subseteq \bigcup_{\bar{j}\sim\bar{i}} Q_{\bar{j}}^{(\lambda)}$ is $[Z(u|\eta)]^{-1}\exp(u\sum_{x\in\mathcal{X}} f(x)\xi_\lambda(x, \mathcal{X}|\eta))$ where $\xi_\lambda(x, \mathcal{X}|\eta)$ stands for the version of the $\lambda$-rescaled packing functional $\xi_\lambda$ in the presence of the marked configuration $\eta$, as discussed above, whereas $Z(u|\eta)$ is the corresponding normalizing constant. The maximum possible number of balls packed in $\bigcup_{\bar{j}\sim\bar{i}} Q_{\bar{j}}^{(\lambda)}$ being of order $O(L^d)$, we have that the logarithm of both this density and $Z(u|\eta)$ is of order $O(uL^d)$ and for the density this holds uniformly in $\mathcal{X}$. Taking $\lambda$ large enough and hence $u$ small enough we make this density fall into the interval $(1 - \varepsilon/2, 1 + \varepsilon/2)$ uniformly in $\mathcal{X}$. Consequently, the probability of any event changing under the considered exponential modification is at most $\varepsilon/2$, as required.

Using the triangle inequality for the total variational norm we conclude that if $u$ is small enough then uniformly over the collection of marked configurations $\eta$, the total variation distance between the law of $\mathcal{P}_{\lambda\kappa}$ restricted to $Q_{\bar{i}}^{(\lambda)}$ and the conditional law of $\mathcal{P}_{\lambda\kappa}^{uf\circ\xi}$ on $Q_{\bar{i}}^{(\lambda)}$ given $\eta$ does not exceed $\varepsilon$. With $\varepsilon$ chosen small enough the assertion of the lemma follows now by Theorem 0.0 in [21] which allows us to construct a coupling of $\mathcal{P}_{\lambda\kappa}$ and $\mathcal{P}_{\lambda\kappa}^{uf\circ\xi}$ such that the disagreement field $(\mathbf{1}_{\mathcal{P}_{\lambda\kappa}\cap Q_{\bar{i}}^{(\lambda)} \neq \mathcal{P}_{\lambda\kappa}^{uf\circ\xi}\cap Q_{\bar{i}}^{(\lambda)}})_{\bar{i}}$ is stochastically dominated by an i.i.d. process of sufficiently low density. $\qquad\square$

We now conclude the proof of (6.12) as follows. Consider the empirical measures $\mu_{\lambda\kappa}^{uf\circ\xi}$ generated by $\mathcal{P}_{\lambda\kappa}^{uf\circ\xi}$ as in (6.10). For $\bar{i}_1,\ldots,\bar{i}_k \in \{1,\ldots,m\}^d$ write

$$m_\lambda^*[\bar{i}_1,\ldots,\bar{i}_k] := \mathbb{E}\prod_{l=1}^k \mu_{\lambda\kappa}^{uf\circ\xi}(Q_{\bar{i}_l}^{(\lambda)}).$$

Using the exponential stabilization of $\xi$ we have the following lemma, which closely corresponds to Lemma 5.2 in [4].

**Lemma 6.2 (Exponential clustering).** *For all integers $p,q \geq 1$ there exist positive constants $A_{p,q}$ and $C_{p,q}$ such that for $\lambda$ large enough one has uniformly in $\bar{i}_1,\ldots,\bar{i}_p; \bar{j}_1,\ldots,\bar{j}_q \in \{1,\ldots,m\}^d$,*

$$|m_\lambda^*[\bar{i}_1,\ldots,\bar{i}_p,\bar{j}_1,\ldots,\bar{j}_q] - m_\lambda^*[\bar{i}_1,\ldots,\bar{i}_p]m_\lambda^*[\bar{j}_1,\ldots,\bar{j}_q]| \leq A_{p,q}\exp\left(-C_{p,q}\min_{r,s}|\bar{i}_r - \bar{j}_s|\right). \tag{6.13}$$

Before proving (6.13), let us first show how it implies (6.12). We do this by adapting the approach of Section 5 of [4] as follows. Put $W := [0,1]^d$ and for all $\lambda > 0$ put $W_\lambda := [0,\lambda^{1/d}]^d$. For any disjoint non-empty $S,T \subset \{1,2,3\}$ define the cluster measure $U_\lambda^{S,T}$ on $W_\lambda^S \times W_\lambda^T$ by

$$U_\lambda^{S,T}(A \times B) = M_\lambda^{S\cup T}(A \times B) - M_\lambda^S(A)M_\lambda^T(B)$$

for all $A \subset W_\lambda^S$ and $B \subset W_\lambda^T$. Here, for any $S := \{s_1,\ldots,s_k\} \subset \{1,2,3\}$, $M_\lambda^S$ is the $|S|$th moment measure on $W^S$, that is the measure on $W^S$ whose discrete probability density function is given by $m_\lambda^*[\bar{i}_{s_1},\ldots,\bar{i}_{s_k}]$ on $\prod_{l=1}^k Q_{\bar{i}_{s_l}}^{(\lambda)}$. As in Lemma 5.1 of [4], for any partition (splitting) of $\{1,2,3\}$ into sets $S,T$, the third cumulant measure $c_{\lambda;\eta hf}^3$ conditioned on the $\sigma$-field generated by products of $Q_{\bar{i}}^{(\lambda)}$, denoted as $\hat{c}_{\lambda;\eta hf}^3$, admits the cluster measure decomposition

$$\hat{c}_{\lambda;\eta hf}^3 = \sum_{(S_1,T_1),(S_2,T_2)} \alpha((S_1,T_1),(S_2,T_2))U_\lambda^{S_1,T_1}M_\lambda^{S_2}M_\lambda^{T_2}, \tag{6.14}$$

where the sum ranges over partitions of $\{1,2,3\}$ consisting of pairings $(S_1,T_1),(S_2,T_2)$ such that $S_1,S_2 \subset S$ and $T_1,T_2 \subset T$ and where $\alpha((S_1,T_1),(S_2,T_2))$ are integer valued pre-factors.

To show (6.12) we need to show

$$\lim_{\lambda\to\infty} \alpha_\lambda \lambda^{-3/2} \int_{W_\lambda^3} f(v_1)f(v_2)f(v_3) \, dc_{\lambda;\eta hf}^3(v_1,v_2,v_3) = 0. \tag{6.15}$$

Taking into account that

$$\int_{W_\lambda^3} f(v_1)f(v_2)f(v_3) \, dc_{\lambda;\eta hf}^3(v_1,v_2,v_3) = (1+\mathrm{o}(1))\int_{W_\lambda^3} f(v_1)f(v_2)f(v_3) \, d\hat{c}_{\lambda;\eta hf}^3(v_1,v_2,v_3)$$

by uniform continuity of $f$ on $[0,1]^d$, it is thus enough to show the following variant of (6.15)

$$\lim_{\lambda\to\infty} \alpha_\lambda \lambda^{-3/2} \int_{W_\lambda^3} f(v_1)f(v_2)f(v_3) \, d\hat{c}_{\lambda;\eta hf}^3(v_1,v_2,v_3) = 0. \tag{6.16}$$

Let $\Delta$ denote the diagonal in $W_\lambda^3$, i.e., $\Delta := \{(v_1, v_2, v_3) \in W_\lambda^3 \colon v_1 = v_2 = v_3\}$. Given $v := (v_1, v_2, v_3) \in W_\lambda^3$, let $D(v)$ denote the Euclidean distance in $(\mathbb{R}^d)^3$ between $v$ and $\Delta$. If the distance in $\mathbb{R}^d$ between all two set partitions of $\{v_1, v_2, v_3\}$ is bounded by $\delta$, then $\operatorname{diam}(\{v_1, v_2, v_3\}) < 3\delta$. However $\operatorname{diam}(\{v_1, v_2, v_3\}) > D(v)/2$, showing that $\delta$ cannot be less than $D(v)/6$. In other words, every triple $v := (v_1, v_2, v_3) \in W_\lambda^3$ can be split into $(v_i)_{i \in S}$ and $(v_j)_{j \in T}$, where $S, T$ is some partition of $\{1, 2, 3\}$, such that the Euclidean distance in $\mathbb{R}^d$ between $(v_i)_{i \in S}$ and $(v_j)_{j \in T}$ exceeds $D(v)/6$. We call this an $(S, T)$ splitting.

We may thus partition $W_\lambda^3$ into a finite number of subsets, say $\sigma(\Pi)$, one for each partition $\Pi := (S, T)$ of $\{1, 2, 3\}$, where all triples $(v_1, v_2, v_3)$ in $\sigma(\Pi)$ admit an $(S, T)$ splitting such that the Euclidean distance in $(\mathbb{R}^d)^3$ between the splittings $(v_i)_{i \in S}$ and $(v_j)_{j \in T}$ exceeds $D(v_1, v_2, v_3)/6$.

We thus have

$$\alpha_\lambda \lambda^{-3/2} \int_{W_\lambda^3} f(v_1) f(v_2) f(v_3) \, \mathrm{d}\hat{c}_{\lambda; \eta h f}^3(v_1, v_2, v_3)$$

$$= \alpha_\lambda \lambda^{-3/2} \sum_\Pi \int_{\sigma(\Pi)} f(v_1) f(v_2) f(v_3) \, \mathrm{d}\hat{c}_{\lambda; \eta h f}^3(v_1, v_2, v_3),$$

where the sum ranges over all partitions $\Pi$ of $\{1, 2, 3\}$. Given a fixed $\Pi$, we now use the corresponding third cumulant representation (6.14) for $\hat{c}_{\lambda; \eta h f}^3$. Exactly as in [4], given this cumulant representation, we decompose the cumulant measure $\mathrm{d}\hat{c}_{\lambda; \eta h f}^3(v_1, v_2, v_3)$ into two measures, one supported by the diagonal $\Delta$ in $W_\lambda^3$ and the other not. Off the diagonal the exponential decay as given by Lemma 6.2 gives a bound of $C\alpha_\lambda \lambda^{-1/2}$. On the diagonal we obtain the same bound, showing (6.16) and (6.12) as desired. Save for the proof of Lemma 6.2, this completes the proof of (6.12) for random sequential packing and spatial birth-growth models.

**Proof of Lemma 6.2.** For simplicity we present the proof for $p = q = 1$, the general case being completely analogous. In fact, we need only $p + q \leq 3$ for our current purposes.

The proof uses Lemma 6.1 together with the exponential decay of the size of subcritical percolation clusters. Put $\bar{i} = \bar{i}_1$, $\bar{j} = \bar{j}_1$ and write $D := d(i, j) := \operatorname{dist}(\bar{i}, \bar{j})/5$ with $\operatorname{dist}(\cdot, \cdot)$ standing here for the usual Euclidean distance in $\mathbb{Z}^d$. For $A \subseteq \mathbb{Z}^d$ write $[A] := [A]_\lambda := \bigcup_{\bar{i} \in A} Q_{\bar{i}}^{(\lambda)}$. We use the coupling of Lemma 6.1 and we declare a sub-cube $Q_{\bar{l}}^{(\lambda)}$ *good* if $\pi_{\bar{l}} = 0$ and *bad* otherwise. Thus $\mathcal{P}_{\lambda \kappa}^{uf \circ \xi}$ and $\mathcal{P}_{\lambda \kappa}$ coincide on good sub-cubes.

We let $\Gamma_{\bar{i}}$ be the maximal subset of the ball $B_{\mathbb{Z}^d}(\bar{i}, D) := \{\bar{l} \in \mathbb{Z}^d, \operatorname{dist}(\bar{l}, \bar{i}) < D\}$ containing $\bar{i}$ and such that we have $\pi_{\bar{l}} = 0$ for all $\bar{l} \in \Gamma_{\bar{i}}^*$ with $\Gamma_{\bar{i}}^*$ standing for the "outer boundary" $\{\bar{l} \in \mathbb{Z}^d \setminus \Gamma_{\bar{i}}, \exists \bar{w} \in \Gamma_{\bar{i}} \bar{l} \sim \bar{w}\}$. If no such subset exists, we put $\Gamma_{\bar{i}} := \emptyset, \Gamma_{\bar{i}}^* := \emptyset$. Likewise, we define $\Gamma_{\bar{j}}$ and $\Gamma_{\bar{j}}^*$. Intuitively speaking, if non-empty, $\Gamma_{\bar{i}}$ is the maximal subset of $B_{\mathbb{Z}^d}(\bar{i}, D)$ completely surrounded by a layer $\Gamma_{\bar{i}}^*$ of *good* sub-cubes.

Clearly, the random sets $\Gamma_{\bar{i}}$ and $\Gamma_{\bar{j}}$ are independent by the coupling Lemma 6.1 and because $B_{\mathbb{Z}^d}(\bar{i}, D) \cap B_{\mathbb{Z}^d}(\bar{j}, D) = \emptyset$ in view of the choice of $D$. Moreover, by the exponential stabilization of $\xi$ we know that the restrictions $\mathcal{P}_{\bar{i}}$ and $\mathcal{P}_{\bar{j}}$ of the point process $\mathcal{P}_{\lambda \kappa}$ (and hence also of $\mathcal{P}_{\lambda \kappa}^{uf \circ \xi}$, by the coupling Lemma 6.1) marked with the values of $\xi_\lambda$ to $[\Gamma_{\bar{i}}^*]$ and $[\Gamma_{\bar{j}}^*]$ respectively can be jointly coupled with their respective independent copies $\mathcal{P}_{\bar{i}}'$ and $\mathcal{P}_{\bar{j}}'$ so that the probability of the event $\{\mathcal{P}_{\bar{i}} \neq \mathcal{P}_{\bar{i}}'\} \cup \{\mathcal{P}_{\bar{j}} \neq \mathcal{P}_{\bar{j}}'\}$ does not exceed the probability that the stabilization sphere at some point of $\mathcal{P}_{\lambda \kappa}$ in $[B_{\mathbb{Z}^d}(\bar{i}, D)]$ extends further than a distance $2LD$ from $L\bar{i}$ or the same happens for $\bar{i}$ replaced with $\bar{j}$ (note that $B_{\mathbb{Z}^d}(\bar{i}, 2D) \cap B_{\mathbb{Z}^d}(\bar{j}, 2D) = \emptyset$ and $B_{\mathbb{Z}^d}(L\bar{i}, 2LD) \cap B_{\mathbb{Z}^d}(L\bar{j}, 2LD) = \emptyset$ by the choice of $D$). Since this probability is of order $\exp(-\Omega(D))$ by exponential stabilization, the total variation distance between the joint law $\mathcal{L}(\mathcal{P}_{\bar{i}}, \mathcal{P}_{\bar{j}})$ and the product $\mathcal{L}(\mathcal{P}_{\bar{i}}) \times \mathcal{L}(\mathcal{P}_{\bar{j}})$ is bounded above by $\exp(-\Omega(D))$.

We will next show that this conclusion also holds for the restrictions $\bar{\mathcal{P}}_{\bar{i}}$ and $\bar{\mathcal{P}}_{\bar{j}}$ of the point process $\mathcal{P}_{\lambda \kappa}^{uf \circ \xi}$ marked with the respective values of $\xi_\lambda$ to the whole sets $[\bar{\Gamma}_{\bar{i}} := \Gamma_{\bar{i}}^* \cup \Gamma_{\bar{i}}]$ and $[\bar{\Gamma}_{\bar{j}}]$ respectively. Indeed, this follows from the simple but important observation that the law of the restriction of $\mathcal{P}_{\lambda \kappa}^{uf \circ \xi}$ to $[\Gamma_{\bar{i}}]$ conditioned on external marked configuration $\eta$ in $[\Gamma_{\bar{i}}]^c$ only depends on this configuration through its restriction to $[\Gamma_{\bar{i}}^*]$ and likewise for $\bar{i}$ replaced by $\bar{j}$ (recall that $L > 4$). Consequently, recalling that $\Gamma_{\bar{i}}^*$ and $\Gamma_{\bar{j}}^*$ are independent, using the processes $\mathcal{P}_i'$ and $\mathcal{P}_j'$ constructed above as marked boundary conditions and *filling* independently the marked point configurations in $[\varGamma_{\bar{i}}]$ and $[\varGamma_{\bar{j}}]$ according to the appropriate conditional distributions we get a joint coupling of $\bar{\mathcal{P}}_{\bar{i}}$ and $\bar{\mathcal{P}}_{\bar{j}}$ with their respective independent copies $\bar{\mathcal{P}}_{\bar{i}}'$ and $\bar{\mathcal{P}}_{\bar{j}}'$ so that the probability of the event $\{\bar{\mathcal{P}}_{\bar{i}} \neq \bar{\mathcal{P}}_{\bar{i}}'\} \cup \{\bar{\mathcal{P}}_{\bar{j}} \neq \bar{\mathcal{P}}_{\bar{j}}'\}$ does not exceed $\exp(-\varOmega(D))$. In other words, the total variation distance between the joint law $\mathcal{L}(\bar{\mathcal{P}}_{\bar{i}}, \bar{\mathcal{P}}_{\bar{j}})$ and the product $\mathcal{L}(\bar{\mathcal{P}}_{\bar{i}}) \times \mathcal{L}(\bar{\mathcal{P}}_{\bar{j}})$ is of order $\exp(-\varOmega(D))$.

To proceed, observe that if $\delta$ in Lemma 6.1 is small enough so that it falls below the critical probability for site percolation on $\mathbb{Z}^d$ with neighborhood relation $\sim$, see [17], by the exponential decay of subcritical percolation cluster size, see Sections 5.2 and 6.3 ibidem, with probability at least $1 - \exp(-\varOmega(D))$ there is no path of *bad* boxes joining either $\bar{i}$ to $[B_{\mathbb{Z}^d}(\bar{i}, 2D)]^c$ or $\bar{j}$ to $[B_{\mathbb{Z}^d}(\bar{j}, 2D)]^c$, whence

$$\mathbf{P}[\bar{i} \in \varGamma_{\bar{i}}, \; \bar{j} \in \varGamma_{\bar{j}}] = \mathbf{P}[\bar{\varGamma}_{\bar{i}} \neq \emptyset, \; \bar{\varGamma}_{\bar{j}} \neq \emptyset] \geq 1 - \exp(-\varOmega(D)).$$

Putting this conclusion together with the independent coupling above and observing that the number of balls packed in $Q_{\bar{i}}^{(\lambda)}$ and $Q_{\bar{j}}^{(\lambda)}$ admits a deterministic upper bound, we conclude the desired relation (6.13) for the correlation functions and hence also (6.12) as required. □

### 6.3. Proof of (6.12) for germ-grain and nearest neighbor functionals

Say that $\xi$ has bounded increments if the increment

$$\Delta^\xi(x, \mathcal{X}) := H^\xi(\mathcal{X} \cup x) - H^\xi(\mathcal{X})$$

admits for all $x$ and $\mathcal{X}$ a deterministic bound $|\Delta^\xi| \leq C_\xi$ for some finite constant $C_\xi$.

In the argument below we shall assume that $\xi$ is either a germ-grain functional (Section 4) or a nearest neighbor functional (Section 5). In the case of the latter, we suppress dependence of $\xi$ on $t$. For these functionals $\xi$ has bounded increments and satisfies $V^\xi(\tau) > 0$, [4, 28].

For the proof of Theorems 4.1 and 5.1 we need a separate argument to establish the required relation (6.12). We also assume that $\kappa$ is bounded away from zero.

We claim that it is enough to prove the following auxiliary lemmas, with $u := \eta h$, stating respectively a *Poisson sandwiching property* and a *mixing property* for the exponentially modified process $\mathcal{P}_{\lambda\kappa}^{uf\circ\xi}$. These are the analogs of Lemmas 6.1 and 6.2. Let $\|\cdot\|_{TV}$ denote the total variation distance.

**Lemma 6.3 (Coupling).** *For $f \in \mathcal{C}([0,1]^d)$ there exist intensities $a := a(f)$, $b := b(f)$, $0 < a < b < \infty$ such that for all $u \in [0,1]$ the point process $\mathcal{P}_{\lambda\kappa}^{uf\circ\xi}$ can be coupled on a common probability space with the homogeneous point processes $\mathcal{P}_{\lambda a}$ and $\mathcal{P}_{\lambda b}$ on $[0,1]^d$ so that a.s. $\mathcal{P}_{\lambda a} \subseteq \mathcal{P}_{\lambda\kappa}^{uf\circ\xi} \subseteq \mathcal{P}_{\lambda b}$.*

**Lemma 6.4 (Exponential clustering).** *For $f \in \mathcal{C}([0,1]^d)$ and for $u$ small enough, for all $k \geq 2$ there exists a constant $C > 0$ such that for $r$ large enough we have*

$$\left\| \mathcal{L}([\mathcal{P}_{\lambda\kappa}^{uf\circ\xi}]_{|\bigcup_{i=1}^k B_{\lambda^{-1/d}r/4}(x_i)]}) - \prod_{i=1}^k \mathcal{L}([\mathcal{P}_{\lambda\kappa}^{uf\circ\xi}]_{|B_{\lambda^{-1/d}r/4}(x_i)}) \right\|_{TV} \leq \exp(-Cr)$$

*for all $x_1, \ldots, x_k \in [0,1]^d$ such that $\min_{i \neq j} \operatorname{dist}(x_i, x_j) > \lambda^{-1/d} r$, with $[\mathcal{P}_{\lambda\kappa}^{uf\circ\xi}]_{|B}$ denoting the restriction of the point process $\mathcal{P}_{\lambda\kappa}^{uf\circ\xi}$ to $B \subseteq [0,1]^d$.*

To complete the proof of (6.12) we recall first that germ grain and nearest neighbor functionals $\xi$ exhibit the exponential stabilization property in the sense of Definition 6.1. We refer the reader to Section 4 of [28] for the proof of this statement. If the process $\mathcal{P}_{\lambda\kappa}^{uf\circ\xi}$ was not dependent then we could directly use Lemma 6.3 with $u := \eta h$ there and with $a := a(f)$ and $b := b(f)$ as given there, and repeat the arguments culminating in Lemma 5.3 in [4] for the usual Poisson point process $\mathcal{P}_{\lambda\kappa}$ replaced with $\mathcal{P}_{\lambda\kappa}^{uf\circ\xi}$ to obtain the exponential decay of correlation functions for $\mu_{\lambda\kappa}^{uf\circ\xi}$. However, while the process $\mathcal{P}_{\lambda\kappa}^{uf\circ\xi}$ is dependent,

Lemma 6.4 shows that its dependency decays exponentially fast with distance, implying that the exponential decay of correlation functions for $\mu_{\lambda\kappa}^{uf\circ\xi}$ has a correction term whose order is bounded by that of the total variation distance expression as in Lemma 6.4. Now, using Lemma 6.4 (note that $u \to 0$ as $\lambda \to \infty$), yields the required decay of the offending third-order term as a particular consequence of the cumulant method developed in [4] which only relies on the exponential decay of correlation functions of all orders, regardless of other particular properties of the considered process, see the discussion in Section 6.2 for a more detailed discussion which applies here as well. Thus, to conclude our argument for germ-grain and nearest neighbor functionals, we establish below Lemmas 6.3 and 6.4.

**Proof of Lemma 6.3.** To provide a constructive proof of the stochastic domination stated in Lemma 6.3 we invoke a general *clan-of-ancestors* algorithm for simulating marked point processes, originally due to Férnandez, Ferrari and Garcia [14, 15, 16]; see also Section 4 in [20] and the references therein. A short description of this algorithm, as provided below, is fully specialized for our particular purposes and corresponds to a very rough version of this simulation scheme, ignoring a number of its essential fine features.

Recalling that $C_\xi$ is the deterministic bound on the increment, and assuming that $\kappa$ is bounded away from zero, we put

$$a := \exp(-u\|f\|_\infty C_\xi) \inf_{x \in [0,1]^d} \kappa(x)$$

and

$$b := \exp(u\|f\|_\infty C_\xi) \sup_{x \in [0,1]^d} \kappa(x).$$

We construct a stationary birth and death point process on $[0,1]^d$ evolving in time $t \in \mathbb{R}$ according to the following dynamics:

- A new point is born at $x \in [0,1]^d$ with intensity $\lambda b\, dx\, dt$,
- An existing point in $[0,1]^d$ dies with intensity $dt$ (i.e. the lifetimes of all points are i.i.d. standard exponential).

The resulting stationary process, denoted in the sequel by $(\Theta_t)_t$, is well defined for all $t \in \mathbb{R}$. It is easily seen that $\Theta_t$ coincides in law on $[0,1]^d$ with $\mathcal{P}_{\lambda b}$ for each $t \in \mathbb{R}$. To proceed, we carry out the following *trimming* procedure for the process $(\Theta_t)_{t \in \mathbb{R}}$. We find the first negative time moment $T_\emptyset$ in the past when $\Theta_t$ becomes empty. Formally,

$$T_\emptyset := \sup\{t < 0 \mid \Theta_t = \emptyset\}.$$

Subsequently we scan all points $\{x_{t_1}, \ldots, x_{t_n}\}$ born in $\Theta_t$, $t \in (T_\emptyset, 0)$, ordered by increasing birth times $t_1 < t_2 < \cdots < t_n$ (note that $t_1 = T_\emptyset$ if $n > 0$) and we accept each $x_{t_i}$, attempting to be born, with probability

$$b^{-1} \exp(uf(x_{t_i}) \Delta^{\xi_\lambda}(x_{t_i}, \mathcal{X}_{t_i})), \tag{6.17}$$

with $\mathcal{X}_{t_i}$ denoting the collection of points accepted by time $t_i$ and alive at that time, otherwise the point $x_{t_i}$ is rejected. Note that the acceptance probability specified in (6.17) falls into $[0, 1]$ by the definition of $b$ and by the bounded increment property of $\xi$. This procedure yields a process $(\Gamma_t)_{t \geq T_\emptyset}$ (in fact the definition of $\Gamma_t$ could easily be extended for all $t \in \mathbb{R}$ by looking further into the past for subsequent evolution epochs separated by periods where the process is empty). It turns out, see [14, 15, 16] and Section 4 in [20], that $\Gamma_0$ coincides in law with $\mathcal{P}_{\lambda\kappa}^{uf\circ\xi}$. It remains to construct the required coupling of $\mathcal{P}_{\lambda a}$, $\mathcal{P}_{\lambda b}$ and $\mathcal{P}_{\lambda\kappa}^{uf\circ\xi}$. To this end, note first that identifying $\mathcal{P}_{\lambda b}$ with $\Theta_0$ and $\mathcal{P}_{\lambda\kappa}^{uf\circ\xi}$ with $\Gamma_0$ yields the required inclusion $\mathcal{P}_{\lambda\kappa}^{uf\circ\xi} \subseteq \mathcal{P}_{\lambda b}$ in view of the obvious relationship $\Gamma_0 \subseteq \Theta_0$. To obtain the required coupling for $\mathcal{P}_{\lambda a}$ assume that the acceptance test (6.17) for points $x_{t_i}$ is carried out by attaching to these points i.i.d. random variables $\beta_i$ uniform on $[0, 1]$ and by declaring a point accepted if the value of $\beta_i$ falls below the corresponding acceptance probability. Then, construct a copy of $\mathcal{P}_{\lambda a}$ by repeating the same acceptance-rejection procedure with the acceptance probabilities in (6.17) replaced by $a/b$ and then by identifying $\mathcal{P}_{\lambda a}$ with the configuration of the resulting *trimmed* process at time 0. Observing that the acceptance probability in (6.17) never falls below $a/b$ we conclude that a.s. $\mathcal{P}_{\lambda a} \subseteq \mathcal{P}_{\lambda b}^{uf\circ\xi}$ as required. □

**Proof of Lemma 6.4.** For better readability we provide the proof when $k = 2$ writing $x$ for $x_1$ and $y$ for $x_2$. The proof for larger $k$ goes exactly along the same lines. The argument below relies on the same *clan-of-ancestors* graphical construction as in the proof of Lemma 6.3. Assume that $u$ is small enough so that

$$b' - a' < \varepsilon$$

for some $\varepsilon > 0$ small enough, with

$$a' := \exp(-uC^\xi \|f\|_\infty) \quad \text{and} \quad b' := \exp(uC^\xi \|f\|_\infty).$$

For formal convenience we replace the birth proposal-acceptance mechanism given above by the following equivalent procedure:

- All points have their lifetimes i.i.d. standard exponential;
- There are two kinds of birth attempts:
  - *regular* birth attempts, which happen with intensity $a'\kappa(x)\,\mathrm{d}x\,\mathrm{d}t$ and are *always* accepted,
  - *exceptional* birth attempts, which happen with intensity $(b' - a')\kappa(x)\,\mathrm{d}x\,\mathrm{d}t$ and are accepted with probability $(b' - a')^{-1}[\exp(uf(x)\Delta^\xi(x, \mathcal{X}_t)) - a']\,\mathrm{d}t$, which clearly lies in $[0, 1]$ by the definition of $a'$ and $b'$. Recall that $\mathcal{X}_t$ stands here for the accepted point configuration at time $t$ when $x$ attempts to be born.

Denoting the resulting stationary process by $(\hat{\Gamma}_t)_{t\geq 0}$ and using again the general setting of [14, 15, 16] and Section 4 of [20] we see that $\hat{\Gamma}_0$ coincides in law with $\mathcal{P}_{\lambda\kappa}^{uf\circ\xi}$ in full analogy with the case of $\Gamma_0$. To proceed, we follow the ideas developed ibidem, constructing an oriented graph on space-time instances of points arising in the above graphical construction, by connecting $y$ to $x$ iff:

- $x$ was created in an *exceptional* birth event,
- $y$ was present at the time where $x$ was born,
- $\mathrm{dist}(x, y) \leq R_{a,b}^{\xi_\lambda}(x, \mathcal{X}_t)$ with $\mathcal{X}_t$ standing for the point configuration present at the moment $t$ of $x$'s birth.

Here the radius of stabilization $R_{a,b}^{\xi_\lambda}(x, \mathcal{X}_t)$ is such that the value of $\xi(x; \mathcal{X}_t)$ remains unaffected by changes outside $B_{R_{a,b}^{\xi_\lambda}(x, \mathcal{X}_t)}(x)$. While the relationship between $R_{a,b}^{\xi_\lambda}(x, \mathcal{X}_t)$ and the stabilization radius in Definition 6.1 may be unclear at the moment, under the coupling constructed in Poisson sandwiching Lemma 6.3 the former will turn out to coincide with the latter. We see that $y$ is connected to $x$ if $y$ could possibly have affected the acceptance status of $x$. We denote by $\mathcal{A}_0[x]$ the union of all directed chains in this graph reaching a given point $x$ alive at time 0 in $\hat{\Gamma}$ and, more generally, for $B \subseteq [0, 1]^d$ we write $\mathcal{A}_0[B]$ for the union of all $\mathcal{A}_0[x]$ with $x \in B$ alive at time 0. Following [14, 15] and [16] these clusters are referred to as *clans of ancestors*, respectively of $x$ and $B$. We see that, for $x, y$ with $\mathrm{dist}(x, y) \geq \lambda^{-1/d}r$, the point process $\mathcal{P}^1 := [\mathcal{P}_{\lambda\kappa}^{uf\circ\xi}]_{|[B_{\lambda^{-1/d}r/16}(x) \cup B_{\lambda^{-1/d}r/16}(y)]}$ can be coupled with independent copies of $\mathcal{P}^2 := [\mathcal{P}_{\lambda\kappa}^{uf\circ\xi}]_{|B_{\lambda^{-1/d}r/16}(x)}$ and $\mathcal{P}^3 := [\mathcal{P}_{\lambda\kappa}^{uf\circ\xi}]_{B_{\lambda^{-1/d}r/16}(y)}$ so that

$$\mathbf{P}[\mathcal{P}^1 \neq \mathcal{P}^2 \cup \mathcal{P}^3] \leq \mathbf{P}\left[\max(\mathrm{diam}\,\mathcal{A}_0[B_{\lambda^{-1/d}r/16}(x)], \mathrm{diam}\,\mathcal{A}_0[B_{\lambda^{-1/d}r/16}(y)]) \geq \lambda^{-1/d}\frac{r}{2}\right]. \quad (6.18)$$

Indeed, this comes from the fact that, by our graphical construction above, the point process $\mathcal{P}_1$ coincides in distribution with the restriction of $\hat{\Gamma}_0$ to $B_{\lambda^{-1/d}r/16}(x) \cup B_{\lambda^{-1/d}r/16}(y)$. Now, again by the construction, conditionally on the event $\{\max(\mathrm{diam}\,\mathcal{A}_0[B_{\lambda^{-1/d}r/16}(x)], \mathrm{diam}\,\mathcal{A}_0[B_{\lambda^{-1/d}r/16}(y)]) \leq \lambda^{-1/d}r/2\}$ the respective restrictions of $\hat{\Gamma}_0$ to $B_{\lambda^{-1/d}r/16}(x)$ and $B_{\lambda^{-1/d}r/16}(y)$ are independent as determined by the realizations of the underlying birth-and-death process with negative time coordinates and with spatial locations restricted respectively to $B_{\lambda^{-1/d}r/2}(x)$ and $B_{\lambda^{-1/d}r/2}(y)$ which are clearly disjoint. This allows us to construct the required coupling.

In the case of germ-grain functionals, this puts us in a position to use a domination by branching processes technique, see ibidem for details, which yields the exponential decay for the tails of the maximum on the right-hand side in (6.18) since the exceptional birth controlling parameter $\varepsilon$ can be made arbitrarily small. This gives Lemma 6.4 for germ-grain functionals.

A slightly different argument is needed for nearest neighbor graphs, where the stabilization radius may admit arbitrarily large values. However, if $R_{a,b}^\xi(x, \mathcal{X}_t) > r$ for a $k$-nearest neighbor graph functional, it clearly means that $B_r(x) \cap \mathcal{X}$ contains at most $k$ points. Consequently, on the event $\{\max(\operatorname{diam}\mathcal{A}_0[B_{\lambda^{-1/d}r/4}(x)], \operatorname{diam}\mathcal{A}_0[B_{\lambda^{-1/d}r/4}]) \geq \lambda^{-1/d}r/2\}$ we must have a connected chain of time–space cylinders, hitting at time 0 either $\partial B_{\lambda^{-1/d}r/4}(x)$ or $\partial B_{\lambda^{-1/d}r/4}(y)$, of length at least $\frac{3}{8}\lambda^{-1/d}r$ and with each cylinder containing less than $k$ points of $(\hat{\Gamma}_t)_{t\geq 0}$. However, using the Poisson sandwiching Lemma 6.3 and its proof we can bound above the probability of existence of such a chain by $\exp(-\Omega(r))$ by the usual discretization and configuration counting argument. Indeed, this is done by splitting the time–space considered in the above graphical construction into fixed-size cubes, using the Poisson sandwiching property to show that the probability of a given collection $\mathcal{C}$ of such cubes containing each less than $k$ points of $(\hat{\Gamma}_t)_{t\geq 0}$ decays like $\exp(-c\operatorname{card}\mathcal{C})$ where $c > 0$ can be made arbitrarily large by suitably adjusting the size of cubes. This allows us to control the parts of the aforementioned chain consisting of large cylinders, while the parts composed of numerous small cylinders are very unlikely to be long if $\varepsilon$ is small enough, by the usual branching process domination technique. We omit the simple but technical details of this argument. Consequently, we get

$$\mathbf{P}\left[\max_{x\in\mathcal{P}^1\cup\mathcal{P}^2\cup\mathcal{P}^3} R_\lambda^\xi[x, \mathcal{P}_{\lambda\kappa}^{uf\xi}] \geq \frac{r}{2}\right] \leq \exp(-Cr), \quad C > 0,$$

which completes the proof of Lemma 6.4. ☐

### 6.4. *Proof of Theorems 2.2 and 2.3*

We prove Theorems 2.2 and 2.3 in the context of the packing measures. The proofs go through without change for the measures described in Theorems 3.1, 4.1 and 5.1.

**Proof of Theorem 2.2.** We will apply the theorems of Gärtner and Ellis combined with the theorem of Dawson and Gärtner dealing with large deviations for projective limits, see Sections 4.5.3 and 4.6 in [10]. Actually we will apply Corollary 4.6.11, part (a), for the family $(\alpha_\lambda^{-1}\lambda^{-1/2}\hat{\mu}_{\lambda\kappa}^\xi)_\lambda$, taking values in the topological vector space $\mathcal{M}([0,1]^d)$. For each $f \in \mathcal{C}([0,1]^d)$ we define $p_f(\nu) := \langle \nu, f \rangle : \mathcal{C}([0,1]^d)' \to \mathbb{R}$, where $\mathcal{C}([0,1]^d)'$ denotes the algebraic dual of $\mathcal{C}([0,1]^d)$ and $\langle \cdot, \cdot \rangle$ denotes the duality brackets between these spaces. Taking the $\sigma$-algebra $\mathcal{B}$ and the maps $p_f$, $f \in \mathcal{C}([0,1]^d)$, we note that $\mathcal{C}([0,1]^d)'$ satisfies Assumption 4.6.8 in [10]. Given the logarithmic moment generating function $\Lambda_{\lambda\kappa,\alpha_\lambda}^\xi(f)$ as in (6.1), we know by Proposition 6.1 that the limit

$$\Lambda_\kappa^\xi(f) := \lim_{\lambda\to\infty} \alpha_\lambda^{-2}\Lambda_{\lambda\kappa,\alpha_\lambda}^\xi(f)$$

exists as an extended real number, and moreover,

$$\Lambda_\kappa^\xi(f) = \frac{1}{2}\int_{[0,1]^d} f^2(x)V^\xi(\kappa(x))\kappa(x)\,dx, \quad f \in \mathcal{C}([0,1]^d).$$

Next we have to show that $\Lambda_\kappa^\xi(\sum_{i=1}^l t_i f_i) : \mathbb{R}^l \to (-\infty, +\infty]$ with $f_1, \ldots, f_l \in \mathcal{C}([0,1]^d)$ and $t_1, \ldots, t_l \in \mathbb{R}$ is, for any $l \in \mathbb{N}$, essentially smooth (for a precise definition see [10], Definition 2.3.5), lower semi-continuous and finite in some neighborhood of 0 (see also Corollary 4.6.14 in [10]). Since $\Lambda_\kappa^\xi$ is everywhere finite, it suffices to show that for every $f, g \in \mathcal{C}([0,1]^d)$, the function $\Lambda_\kappa^\xi(f + tg) : \mathbb{R} \to \mathbb{R}$ is differentiable at $t = 0$, implying that $\Lambda_\kappa^\xi(\cdot)$ is Gateaux differentiable. Using the boundedness of $f$ and $g$ we may interchange the order of differentiation and integration to obtain

$$\frac{\mathrm{d}}{\mathrm{d}t}\Lambda_\kappa^\xi(f+tg)|_{t=0} = \int_{[0,1]^d} f(x)g(x)V^\xi(\kappa(x))\kappa(x)\,\mathrm{d}x,$$

which is well defined. Consequently, all the conditions of part (a) of Corollary 4.6.11 in [10] are satisfied, and hence the family $(\alpha_\lambda^{-1}\lambda^{-1/2}\tilde{\mu}_{\lambda\kappa}^\xi)_\lambda$ satisfies the LDP in $\mathcal{C}([0,1]^d)'$ with respect to the weak topology, with speed $\alpha_\lambda^2$ and a convex, good rate function

$$\Lambda_{\xi,\kappa}^*(\nu) := \sup_{f\in\mathcal{C}([0,1]^d)} \{\langle f,\nu\rangle - \Lambda_\kappa^\xi(f)\} \quad \text{for } \nu\in\mathcal{C}([0,1]^d)'.$$

Now we restrict the LDP. We will show that $\Lambda_{\xi,\kappa}^*(\nu) < \infty$ implies that $\nu$ is a continuous linear form. Assume that $\Lambda_{\xi,\kappa}^*(\nu) < \infty$. Then for all $f\in\mathcal{C}([0,1]^d)$, $f\neq 0$, we obtain

$$\left\langle \frac{f}{\|f\|},\nu\right\rangle \leq \Lambda_{\xi,\kappa}^*(\nu) + \Lambda_\kappa^\xi\left(\frac{f}{\|f\|}\right) \leq \Lambda_{\xi,\kappa}^*(\nu) + \Lambda_\kappa^\xi(1).$$

Now by assumption the right-hand side is finite and hence $\langle\nu,f\rangle \leq K\|f\|$. Considering $-f$, we get $|\langle\nu,f\rangle| \leq K\|f\|$. Thus $\nu$ is a continuous linear form. Since $[0,1]^d$ is compact, we can apply Riesz's representation theorem which implies that $\nu$ can be represented as a $\mathbb{R}$-valued measure on $[0,1]^d$, e.g. $\nu\in\mathcal{M}([0,1]^d)$. We can apply Lemma 4.1.5 in [10] to obtain the LDP in the space $\mathcal{M}([0,1]^d)$, endowed with the weak topology. The representation of the rate function is proved in Lemma 6.5, hence the proof of Theorem 2.2 is complete. $\square$

**Lemma 6.5.** *The identity*

$$\Lambda_{\xi,\kappa}^*(\cdot) = I_\kappa^\xi(\cdot)$$

*holds over $\mathcal{M}([0,1]^d)$. Moreover, we obtain*

$$\Lambda_{\xi,\kappa}^*(\nu) = \sup_{f\in B([0,1]^d)} \{\langle f,\nu\rangle - \Lambda_\kappa^\xi(f)\}, \quad \nu\in\mathcal{M}([0,1]^d), \tag{6.19}$$

*where $B([0,1]^d)$ denotes the collection of bounded measurable functions $f\colon[0,1]^d\to\mathbb{R}$.*

**Proof.** Let us define the measure $\mu(\xi,\kappa)$ by $\mathrm{d}\mu(\xi,\kappa) := V^\xi(\kappa(x))\kappa(x)\,\mathrm{d}x$. For a fixed $\nu\in\mathcal{M}([0,1]^d)$ and $f\in B([0,1]^d)$ there exists a sequence $f_n\in\mathcal{C}([0,1]^d)$, $n\in\mathbb{N}$, such that $\lim f_n = f$ both in $L_1(\mu(\xi,\kappa))$ and in $L_1(\nu)$ (truncating each $f_n$ to the bounded range of $f$). Consequently, there exists a sequence $(f_n)_n\subset\mathcal{C}([0,1]^d)$ such that

$$\lim_{n\to\infty}\left(\int_{[0,1]^d} f_n\,\mathrm{d}\nu - \Lambda_\kappa^\xi(f_n)\right) = \int_{[0,1]^d} f\,\mathrm{d}\nu - \Lambda_\kappa^\xi(f).$$

Since $\nu\in\mathcal{M}([0,1]^d)$ and $f\in B([0,1]^d)$ are arbitrary, the definition of $\Lambda_{\xi,\kappa}^*$ agrees with (6.19) over $\mathcal{M}([0,1]^d)$.

To proceed, let us assume that $\nu\in\mathcal{M}([0,1]^d)$ is chosen such that $I_\kappa^\xi(\nu) < \infty$. This is true by definition only for those $\nu$ with a density $\varrho$ with respect to the measure $\mu(\xi,\kappa)$. Hence $I_\kappa^\xi(\nu) = \frac{1}{2}\int\varrho^2\,\mathrm{d}\mu(\xi,\kappa)$. Since for every $f\in\mathcal{C}([0,1]^d)$

$$f\varrho \leq \frac{1}{2}f^2 + \frac{1}{2}\varrho^2,$$

we obtain

$$\int f\,\mathrm{d}\nu - \frac{1}{2}\int f^2\,\mathrm{d}\mu(\xi,\kappa) = \int f\varrho\,\mathrm{d}\mu(\xi,\kappa) - \frac{1}{2}\int f^2\,\mathrm{d}\mu(\xi,\kappa) \leq \frac{1}{2}\int\varrho^2\,\mathrm{d}\mu(\xi,\kappa) = I_\kappa^\xi(\nu)$$

and therefore $\Lambda_{\xi,\kappa}^*(\cdot) \leq I_\kappa^\xi(\cdot)$ over the whole $\mathcal{M}([0,1]^d)$. To get the converse inequality take $\nu \in \mathcal{M}([0,1]^d)$ such that

$$\mathrm{d}\nu := f\,\mathrm{d}\mu(\xi,\kappa), \quad f \in L_1(\mu(\xi,\kappa)),$$

and define a sequence $f_n \in B([0,1]^d)$, $n \in \mathbb{N}$, by $f_n(x) := \mathrm{sign}(f(x))\min(|f(x)|,n)$. Then we simply obtain the identities

$$\begin{aligned}
I_\kappa^\xi(\nu) &= \int_{[0,1]^d} \left(\frac{\mathrm{d}\nu}{\mathrm{d}\mu(\xi,\kappa)}\right)^2 \mathrm{d}\mu(\xi,\kappa) - \frac{1}{2}\int_{[0,1]^d} \left(\frac{\mathrm{d}\nu}{\mathrm{d}\mu(\xi,\kappa)}\right)^2 \mathrm{d}\mu(\xi,\kappa) \\
&= \int_{[0,1]^d} \frac{\mathrm{d}\nu}{\mathrm{d}\mu(\xi,\kappa)}\,\mathrm{d}\nu - \frac{1}{2}\int_{[0,1]^d} \left(\frac{\mathrm{d}\nu}{\mathrm{d}\mu(\xi,\kappa)}\right)^2 \mathrm{d}\mu(\xi,\kappa) \\
&= \int_{[0,1]^d} f\,\mathrm{d}\nu - \frac{1}{2}\int_{[0,1]^d} f^2\,\mathrm{d}\mu(\xi,\kappa) \\
&= \lim_{n\to\infty} \int_{[0,1]^d} f_n\,d\nu - \frac{1}{2}\int_{[0,1]^d} f_n^2\,\mathrm{d}\mu(\xi,\kappa) \leq \Lambda_{\xi,\kappa}^*(\nu).
\end{aligned}$$

To complete the proof of Lemma 6.5 it remains to show that $\Lambda_{\xi,\kappa}^*(\nu) = +\infty$ for $\nu$ which are not absolutely continuous with respect to $\mu(\xi,\kappa)$. But this follows immediately from the fact that for such $\nu$ we can find $f \in B([0,1]^d)$ with $\langle f, \nu \rangle$ arbitrarily large and $\Lambda_\kappa^\xi(f) = \frac{1}{2}\int_{[0,1]^d} f^2\,\mathrm{d}\mu(\xi,\kappa)$ arbitrarily small. This finishes the proof of Lemma 6.5. $\qquad\square$

**Proof of Theorem 2.3.** Denote by $p_{f_1,\dots,f_l}\colon \mathcal{C}([0,1]^d)' \to \mathbb{R}^l$ the function

$$p_{f_1,\dots,f_l}(\nu) := (\langle f_1, \nu \rangle, \dots, \langle f_l, \nu \rangle).$$

Note that the limiting logarithmic moment generating function associated with the family

$$(\alpha_\lambda^{-1}\lambda^{-1/2}\bar{\mu}_{\lambda\kappa}^\xi \circ p_{f_1,\dots,f_l}^{-1})_\lambda$$

is the function

$$h(t) := \Lambda_\kappa^\xi\left(\sum_{i=1}^l t_i f_i\right)\colon \mathbb{R}^l \to (-\infty, +\infty], \quad t = (t_1,\dots,t_l).$$

But we have checked in the proof of Theorem 2.2, that $h$ is essentially smooth, lower semi-continuous and finite in some neighborhood of the origin. Hence the Gärtner–Ellis theorem [10], Theorem 2.3.6, implies that these measures satisfy a LDP in $\mathbb{R}^l$ with speed $\alpha_\lambda^2$ and with the good rate function $I_{f_1,\dots,f_l}\colon \mathbb{R}^l \to [0,\infty]$, where for any $\nu \in \mathcal{C}_b([0,1]^d)'$, putting $s_i := \langle f_i, \nu \rangle$ one gets

$$I_{f_1,\dots,f_l}^\xi((s_1,\dots,s_l)) = \sup_{t_1,\dots,t_l \in \mathbb{R}} \left(\sum_{i=1}^l t_i s_i - \Lambda_\kappa^\xi\left(\sum_{i=1}^l t_i f_i\right)\right). \tag{6.20}$$

It easy to see that

$$\Lambda_\kappa^\xi\left(\sum_{i=1}^l t_i f_i\right) = \Lambda_\kappa^\xi(\langle t, (f_1,\dots,f_l)\rangle) = \langle t, C(\xi,\kappa,f_1,\dots,f_l)t\rangle.$$

It follows that $I_{f_1,\dots,f_l}$ in (6.20) coincides with $I_{\kappa,f_1,\dots,f_l}^\xi$ in (2.4). This completes the proof. $\qquad\square$

## 7. Proof of the LIL (Theorem 2.4)

We prove Theorem 2.4 for the packing measures. The proofs go through without change for the measures described in Theorems 3.1, 4.1 and 5.1.

Recall that in the argument below we consider $\mathcal{M}([0,1]^d)$ endowed with the topology metrized by $\mathrm{dist}_W(\cdot,\cdot)$ as given in (2.7). The first step of the proof of Theorem 2.4 makes standard use of Theorem 2.2 in order to show that $(\zeta_{\kappa\lambda}^\xi)_\lambda$ is almost surely compact and its set of accumulation points is almost surely contained in $\mathcal{K}_\kappa^\xi$. Next, in the second part of the proof we use the rapid decay of dependencies, as concluded from the exponential stabilization in Lemma 7.2 below, to construct an appropriate coupling so that almost surely each point of $\mathcal{K}_\kappa^\xi$ is attained as an accumulation point of $(\zeta_{\lambda\kappa}^\xi)_\lambda$.

### 7.1. Accumulation points

Recall that $\mathcal{K}_\kappa^\xi$ is compact since $I_\kappa^\xi$ is a good rate function. Thus, to establish the relative compactness of $(\zeta_{\lambda\kappa}^\xi)_\lambda$ and to prove that almost surely all accumulation points of $(\zeta_{\lambda\kappa}^\xi)_\lambda$ fall into $\mathcal{K}_\kappa^\xi$ is enough to show that almost surely

$$\limsup_{\lambda\to\infty}\mathrm{dist}_W(\zeta_{\lambda\kappa}^\xi,\mathcal{K}_\kappa^\xi)=0. \tag{7.1}$$

We use the following lemma, which is a straightforward modification of Lemma 1.4.3 in [11].

**Lemma 7.1.** *To establish* (7.1) *it is enough to show that for each $s>1$ we have almost surely*

$$\limsup_{k\to\infty}\mathrm{dist}_W(\zeta_{s^k\kappa}^\xi,\mathcal{K}_\kappa^\xi)=0. \tag{7.2}$$

To proceed with the proof of (7.1), fix $s>1$ and choose arbitrary $\eta>0$. Clearly,

$$\inf\{I_\kappa^\xi(\theta)\mid\mathrm{dist}_W(\theta,\mathcal{K}_\kappa^\xi)\geq\eta\}>1.$$

In particular, in view of the moderate deviation principle upper bound in Theorem 2.2 we have for $k$ large enough,

$$\mathbf{P}[\mathrm{dist}_W(\zeta_{s^k\kappa}^\xi,\mathcal{K}_\kappa^\xi)\geq\eta]\leq\exp(-\alpha_{s^k}^2(1+\delta))=(k\log s)^{-(1+\delta)} \tag{7.3}$$

with some $\delta>0$. Hence, in view of the Borel–Cantelli lemma, the event

$$\{\mathrm{dist}_W(\zeta_{s^k\kappa}^\xi,\mathcal{K}_\kappa^\xi)\geq\eta\}$$

occurs almost surely at most a finite number of times. Consequently, almost surely each accumulation point $\theta^*$ of $(\zeta_{s^k\kappa}^\xi)_{k=1}^\infty$ in $\mathcal{M}([0,1]^d)$ has to satisfy $\mathrm{dist}_W(\theta^*,\mathcal{K}_\kappa^\xi)<\eta$. As $\eta$ was arbitrary, we conclude (7.2) and hence (7.1).

### 7.2. Conclusion of the proof of the LIL

The following lemma, stating rapid decay of dependencies (exponential $\alpha$-mixing in fact) between $\mu_{\lambda\kappa}^\xi(A)$ and $\mu_{\lambda\kappa}^\xi(B)$ for separated $A,B\subseteq[0,1]^d$, will be an important tool in the sequel of our argument.

**Lemma 7.2.** *Let $A,B\subseteq[0,1]^d$ be measurable sets satisfying $\mathrm{dist}(A,B)>r$ and denote by $\mathcal{F}_A^\lambda$ and $\mathcal{F}_B^\lambda$ the sigma fields generated by the restrictions of the random measure $\mu_{\lambda\kappa}^\xi$ to $A$ and $B$, respectively. Then, for any two events $\mathcal{E}_1\in\mathcal{F}_A^\lambda$ and $\mathcal{E}_2\in\mathcal{F}_B^\lambda$ we have*

$$|\mathbf{P}[\mathcal{E}_1\cap\mathcal{E}_2]-\mathbf{P}[\mathcal{E}_1]\mathbf{P}[\mathcal{E}_2]|=\mathrm{O}\left(\lambda L\int_{A\cup B}\kappa(x)\,\mathrm{d}x\exp\left(-\alpha\lambda^{1/d}\frac{r}{2}\right)\right).$$

**Proof.** Write

$$\mathcal{O}_A := \left\{ x \in A \cap \mathcal{P}_{\lambda\kappa} \colon R^\xi_{\lambda;\kappa}(x) \leq \lambda^{1/d}\frac{r}{2} \right\}$$

and let $\mathcal{O}_B$ be defined analogously. Since $\mathrm{dist}(A,B) > r$, it follows by the definition of the re-scaled measure $\mu^\xi_{\lambda\kappa}$ and by stabilization that

$$\mathbf{P}[\mathcal{E}_1 \cap \mathcal{E}_2 | \mathcal{O}_A \cap \mathcal{O}_B] = \mathbf{P}[\mathcal{E}_1 | \mathcal{O}_A \cap \mathcal{O}_B]\mathbf{P}[\mathcal{E}_2 | \mathcal{O}_A \cap \mathcal{O}_B], \tag{7.4}$$

i.e. $\mathcal{E}_1$ and $\mathcal{E}_2$ are conditionally independent on $\mathcal{O}_A \cap \mathcal{O}_B$. Moreover, by (6.3) we get

$$\mathbf{P}[\mathcal{O}_A \cup \mathcal{O}_B] \geq 1 - \mathbb{E}\left[\sum_{x \in A \cap \mathcal{P}_{\lambda\kappa}} \mathbf{1}_{\{R^\xi_{\lambda;\kappa}(x) > \lambda^{1/d}r\}}\right] - \mathbb{E}\left[\sum_{x \in B \cap \mathcal{P}_{\lambda\kappa}} \mathbf{1}_{\{R^\xi_{\lambda;\kappa}(x) > \lambda^{1/d}r\}}\right]$$

$$\geq 1 - \lambda L \int_{A \cup B} \kappa(x)\,\mathrm{d}x \exp\left(-\alpha\lambda^{1/d}\frac{r}{2}\right).$$

Consequently, in view of (7.4),

$$|\mathbf{P}[\mathcal{E}_1 \cap \mathcal{E}_2] - \mathbf{P}[\mathcal{E}_1]\mathbf{P}[\mathcal{E}_2]|$$

$$\leq |\mathbf{P}[\mathcal{E}_1 \cap \mathcal{E}_2] - \mathbf{P}[\mathcal{E}_1 \cap \mathcal{E}_2 | \mathcal{O}_A \cap \mathcal{O}_B]| + |\mathbf{P}[\mathcal{E}_1]\mathbf{P}[\mathcal{E}_2] - \mathbf{P}[\mathcal{E}_1 | \mathcal{O}_A \cap \mathcal{O}_B]\mathbf{P}[\mathcal{E}_2 | \mathcal{O}_A \cap \mathcal{O}_B]|$$

$$= \mathrm{O}\left(\lambda L \int_{A \cup B} \kappa(x)\,\mathrm{d}x \exp\left(-\alpha\lambda^{1/d}\frac{r}{2}\right)\right)$$

as required. $\qquad\square$

To proceed, consider the following coupling of the Poisson point processes $\mathcal{P}_{\lambda\kappa}$ on $(\Omega, \mathcal{F}, \mathbf{P})$. Let $\Pi_1$ be a homogeneous Poisson point process on $(\mathbb{R}^+)^d \times \mathbb{R}^+$ with intensity 1. Then $\mathcal{P}_{\lambda\kappa}$ can be identified with the point process

$$\{\lambda^{-1/d}x \colon \exists t > 0\,(x,t) \in \Pi_1, x \in [0, \lambda^{1/d}]^d, t \leq \kappa(\lambda^{-1/d}x)\}. \tag{7.5}$$

It is worth noting that for constant $\kappa(x) \equiv \kappa$ the above coupling corresponds to observing always the same homogeneous Poisson point process of intensity $\kappa$ on $(\mathbb{R}^+)^d$, while successively increasing the observation window $[0, \lambda^{1/d}]^d$ and then re-scaling the observations onto $[0,1]^d$, thus getting copies of $\mathcal{P}_{\lambda\kappa}$. Clearly, this is a natural coupling construction appearing in a number of applications.

We shall show that with the coupling (7.5), almost surely each $\theta \in \mathcal{K}^\xi_\kappa$ is attained as an accumulation point of $\zeta^\xi_{\lambda\kappa}$ for $\lambda \to \infty$, in other words, almost surely for each $\theta \in \mathcal{K}^\xi_\kappa$ there exists a sequence $\lambda^{(\theta)}_k \to \infty$ with $\zeta^\xi_{\lambda^{(\theta)}_k \kappa} \to \theta$. Clearly, in view of the separability of $\mathcal{K}^\xi_\kappa$, it is enough to prove, separately for each $\theta \in \mathcal{K}^\xi_\kappa$, that, almost surely,

$$\liminf_{\lambda \to \infty} \mathrm{dist}_W(\zeta^\xi_{\lambda\kappa}, \theta) = 0. \tag{7.6}$$

To this end, fix $\theta \in \mathcal{K}^\xi_{\lambda\kappa}$, choose arbitrary $\eta > 0$ and let $f \in \mathcal{C}([0,1]^d)$ be a function with its support contained in $[\delta, 1-\delta]^d$ for some $\delta > 0$. Without loss of generality we assume that $\delta < 1/2$ and we set

$$\rho := \frac{1-\delta}{\delta} > 1.$$

We claim that in order to establish (7.6) it is enough to show that, almost surely,

$$\liminf_{k \to \infty} |\langle f, \zeta^\xi_{\rho^{kd}\kappa}\rangle - \langle f, \theta\rangle| = 0. \tag{7.7}$$

Indeed, since $f$ was arbitrary, by a standard diagonal argument (7.7) yields almost sure existence of a (random) subsequence $k_n$ with

$$\lim_{k_n \to \infty} \langle f_m, \zeta^\xi_{\rho^{k_n}{}_\kappa} \rangle = \langle f_m, \theta \rangle \tag{7.8}$$

for all $f_m$ running through a countable collection of continuous functions with compact supports bounded away from $\partial[0,1]^d$, uniformly dense in the set of all continuous functions on $[0,1]^d$ assuming the value 0 at the boundary $\partial[0,1]^d$. Since $(\zeta^\xi_{\lambda\kappa})_\lambda$ is almost surely relatively compact in $(\mathcal{M}([0,1]^d), \mathrm{dist}_W)$ as argued in Section 7.1 above, setting $W := W \cup \{f_m\}_{m=1}^\infty$ for $W$ in (2.7) and using (7.8) (note that enlarging $W$ while keeping it countable brings no loss of generality and can only strengthen our result as refining the induced topology) we conclude that almost surely there exists a random subsequence $k'_n$ with $\zeta^\xi_{\rho^{k'_n}{}_\kappa}$ converging weakly to some (random) measure $\theta'$ satisfying almost surely

$$\langle f_m, \theta \rangle = \langle f_m, \theta' \rangle$$

for all $f_m$ as above. Recalling that $I^\xi_\kappa(\theta') \le 1$ and, consequently, $\theta'(\partial[0,1]^d) = 0$, we conclude that $\theta' = \theta$ almost surely, which shows that $\theta$ arises almost surely as an accumulation point of $(\zeta^\xi_{\rho^k \kappa})_k$, as required.

It remains to establish (7.7). To this end, consider the sequence of disjoint cubes

$$Q_k := [\delta\rho^k, (1-\delta)\rho^k]^d$$

and, fixing some $\eta > 0$, denote by $E_k := E_k^{(f,\eta)}$ the event

$$E_k := \{ |\langle f, \zeta^\xi_{\rho^{kd}\kappa} \rangle - \langle f, \theta \rangle| \ge \eta \}.$$

Applying the moderate deviation principle as stated in Theorem 2.2 we see that, for $k$ large enough,

$$\mathbf{P}[E_k^c] = \mathbf{P}[|\langle f, \zeta^\xi_{\rho^{kd}\kappa} \rangle - \langle f, \theta \rangle| < \eta] \ge \exp(-\alpha^2_{\rho^{kd}}[\inf\{I^\xi_\kappa(\gamma) \mid |\langle f, \gamma \rangle - \langle f, \theta \rangle| < \eta\} + \varepsilon]),$$

with $\varepsilon > 0$ chosen so that

$$\beta := \inf\{I^\xi_\kappa(\gamma) \mid |\langle f, \gamma \rangle - \langle f, \theta \rangle| < \eta\} + \varepsilon < 1$$

(note that such a choice is possible as the rate function $I^\xi_\kappa$ has the property that in an arbitrarily small open neighborhood of any $\theta \in \mathcal{K}^\xi_\kappa$ one can always find $\gamma$ with $I^\xi_\kappa(s) < I^\xi_\kappa(\theta) \le 1$). Consequently, we obtain

$$\mathbf{P}[E_k^c] \ge (kd\log\rho)^{-\beta}. \tag{7.9}$$

To proceed, consider the event

$$A_m := \bigcap_{k \ge m} E_k.$$

It is clear that in order to establish (7.7) it is enough to show that

$$\mathbf{P}[A_m] = 0 \tag{7.10}$$

for all $m > 0$. Taking into account that the distance between $Q_k$ and $\bigcup_{j=1}^k Q_j$ is larger than $\delta\rho^{k-1}$, which corresponds to $\lambda^{-1/d}\delta\rho^{k-1}$ with $\lambda = \rho^{kd}$ under the re-scaling in the definition of $\mu^\xi_{\lambda\kappa}$, we can apply Lemma 7.2 and use (7.9) to conclude for all $k > m$ that

$$\mathbf{P}[A_m] \le \mathbf{P}[E_m \cap \cdots \cap E_k] \le \mathbf{P}[E_m \cap \cdots \cap E_{k-1}]\mathbf{P}[E_k] + \mathrm{O}\left(\lambda L \int_{[0,1]^d} \kappa(x)\,\mathrm{d}x \exp\left(-\alpha\lambda^{1/d}\lambda^{-1/d}\delta\rho^{k-1}\frac{r}{2}\right)\right)$$

$$\le \mathbf{P}[E_m \cap \cdots \cap E_k](1 - (kd\log\rho)^{-\beta}) + \mathrm{O}\left(\lambda L \int_{[0,1]^d} \kappa(x)\,\mathrm{d}x \exp\left(-\alpha\delta\rho^{k-1}\frac{r}{2}\right)\right).$$

Consequently, we obtain

$$\mathbf{P}[A_m] \le \exp\left(-\sum_{k=m}^{\infty}(kd\log\rho)^{-\beta}\right) + O\left(\sum_{k=m}^{\infty}\rho^{kd}L\exp\left(-\alpha\delta\rho^{k-1}\frac{r}{2}\right)\int_{[0,1]^d}\kappa(x)\,dx\right). \tag{7.11}$$

Since the first term on the right-hand side of the inequality (7.11) is zero while the rightmost series converges, we conclude that

$$\lim_{m\to\infty}\mathbf{P}[A_m] = 0.$$

Since $A_m$ is an increasing sequence of events, this yields (7.10) and, consequently, also the relations (7.7) and (7.6). The proof of Theorem 2.4 is hence complete.

## 8. Proof of the LIL (binomial case, Theorem 2.6)

We prove Theorem 2.6 for the packing measures (2.10). The proofs go through without change for the corresponding measures described in Theorems 3.1, 4.1 and 5.1. Put

$$\gamma := \int_{[0,1]^d} f(x)\delta^{\xi}(\kappa(x))\kappa(x)\,dx.$$

By coupling arguments and the bounded moments condition (see e.g. (6.10) in [4] or (5.32) in [26]) there exists a coupling of the binomial measures $\rho_n^{\xi}$ and the Poisson measures $\mu_{n\kappa}^{\xi}$, $n = 1, 2, \ldots,$ such that

$$\mathbb{E}\left[\left(n^{-1/2}\langle f, \mu_{n\kappa}^{\xi}\rangle - \langle f, \rho_{n,\kappa}^{\xi}\rangle - \gamma\left(\sum_{j=1}^{n}\eta_j - n\right)\right)^2\right] \to 0, \tag{8.1}$$

where $\eta_1, \eta_2, \ldots,$ is a sequence of i.i.d. Poisson(1) random variables with $(\eta_j)_{j\ge 1}$ independent of $(\rho_{n,\kappa}^{\xi})_{n\ge 1}$. In particular, we get from (8.1) that

$$\mathbb{E}\left(\langle f, \zeta_{n\kappa}^{\xi}\rangle - \langle f, \theta_{n,\kappa}^{\xi}\rangle - \gamma\alpha_n^{-1}n^{-1/2}\left(\sum_{j=1}^{n}\eta_j - n\right)\right)^2 \to 0. \tag{8.2}$$

The classical law of the iterated logarithm yields almost surely

$$\limsup_{n\to\infty}\alpha_n^{-1}n^{-1/2}\left(\sum_{j=1}^{n}\eta_j - n\right) = \sqrt{2} \tag{8.3}$$

and

$$\liminf_{n\to\infty}\alpha_n^{-1}n^{-1/2}\left(\sum_{j=1}^{n}\eta_j - n\right) = -\sqrt{2}. \tag{8.4}$$

Thus, in view of the independence of $(\rho_{n,\kappa}^{\xi})_n$ and $(\eta_j)_j$, putting (8.2), (8.3) and (8.4) together we conclude that a violation with a positive probability of either of the inequalities (2.13) and (2.14) would lead to a violation of (2.8) or (2.9). This completes the proof of Theorem 2.6.

## Acknowledgments

As an alternative to our current approach one could try to use the cumulant-based methods which have been widely developed in the Russian probability literature, see e.g. [32] and the references therein. However, we have chosen the techniques in the proof of Proposition 6.1 as they are presented in a self-contained way and are directly specialized for the particular multidimensional continuum measure-valued setting considered in this paper, as opposed to traditional formulations of the aforementioned results, which would require considerable effort to make them directly applicable to our setting. The third author wishes to express his gratitude to A. V. Nagaev and Z. S. Szewczak for pointing out the above references and for valuable comments. The second and fourth authors greatly appreciate the hospitality of the Institute for Mathematical Sciences at the National University of Singapore, where parts of this work began.